\documentclass[11 pt,reqno]{amsart}
\usepackage{amsmath,amsthm,amsbsy,amssymb}

\usepackage{graphicx}

 \evensidemargin
 \oddsidemargin
\newcommand{\smooth}{\mathcal{S}}
\renewcommand{\leq}{\leqslant}
\renewcommand{\geq}{\geqslant}

 \DeclareMathOperator{\E}{{\mathbb{E}}}
 
 \DeclareMathOperator{\Ci}{Ci}

\newcommand{\R}{{\mathbb{R}}}

\newcommand{\T}{{\mathcal{T}}}
\renewcommand{\Re}{{\rm{Re}}}
\renewcommand{\Im}{{\rm{Im}}}

\renewcommand{\a}{\alpha}

\newcommand{\g}{\gamma}

\renewcommand{\t}{\theta}

\renewcommand{\i}{{\mathrm{i}}}   
\renewcommand{\d}{{\mathrm{d}}}   

\newtheorem{thm}{Theorem}
\newtheorem*{thm*}{Theorem}

\newtheorem*{cor*}{Corollary}
\newtheorem{lem}{Lemma}
\newtheorem{conj}{Conjecture}
\theoremstyle{remark}
\newtheorem*{rem}{Remark}

\begin{document}

\title{A Hybrid Euler-Hadamard product formula for the Riemann zeta function}

\author{S.M. Gonek}
 \address{Department of Mathematics, University of Rochester, Rochester, NY 14627, USA}
 \email{gonek@math.rochester.edu}

\author{C.P. Hughes}
 \address{Department of Mathematics, University of Michigan, Ann Arbor, MI 48109, USA}
 \email{hughes@aimath.org}

\author{J.P. Keating}
\address{School of Mathematics, University of Bristol, University Walk, Bristol,
BS8 1TW, UK}
 \email{j.p.keating@bristol.ac.uk}

\date{November 5, 2005}

\begin{abstract}
We use a smoothed version of the explicit formula to find
an approximation to the Riemann zeta function
as a product over its nontrivial zeros
multiplied by a product over the primes.
We model the first product
by characteristic polynomials of random matrices.
This provides a statistical model of the zeta function
that involves the primes in a natural way. We then employ
the model in a heuristic calculation of the moments of the modulus
of the zeta function on the critical line. This calculation illuminates recent
conjectures for these moments based on connections with random matrix
theory.
\end{abstract}

\maketitle

\section{Introduction}\label{sect:intro}

An important theme in the study of the Riemann zeta function,
$\zeta(s)$, has been
the estimation of the mean values (or moments)
\begin{equation*}
I_k(T) = \frac 1T \int_{0}^{T} \left| \zeta(\tfrac 12 + \i t)
\right|^{2k} \;\d t \; .
\end{equation*}
These have applications to bounding the order of $\zeta(s)$ in the
critical strip as well as to estimating the possible number of zeros of the
zeta function off the critical line. Moreover, the techniques developed in these
problems, in addition to being interesting in their own right,
have been used to estimate mean values of other important functions
in analytic number theory, such as Dirichlet polynomials.

In 1918 Hardy and Littlewood~\cite{HardyLittle} proved that
\begin{equation*}
I_{1}(T) \sim \log T
\end{equation*}
as $T \to \infty$. Eight years later, in 1926, Ingham~\cite{Ing} showed that
\begin{equation*}
I_{2}(T) \sim \frac1{2\pi^{2}}(\log T)^{4} \;.
\end{equation*}
There are no proven asymptotic results for  $I_{k}$ when $ k > 2$, although
it has long been conjectured that
\begin{equation*}
I_{k}(T) \sim c_{k}(\log T)^{k^{2}}
\end{equation*}
for some positive constant $c_{k}$. Conrey and
Ghosh (unpublished) cast this in a more precise form, namely,
\begin{equation*}
I_{k}(T) \sim \frac {a(k)g(k)}{\Gamma(k^{2}+1)}(\log T)^{k^{2}}\;,
\end{equation*}
where
\begin{equation}\label{eq:a(k)}
a(k) = \prod_{p}\left(
\left(1-\frac{1}{p}\right)^{k^2} \sum_{m=0}^\infty
\left(\frac{\Gamma(m+k)}{m! \ \Gamma(k)}\right)^{\!2} p^{-m} \right)\; ,
\end{equation}
the product being taken over all prime numbers, and $g(k)$ is an
integer when $k$ is an integer. The results of Hardy--Littlewood
and Ingham give $g(1)=1$ and $g(2)=2$, respectively. However,
until recently, no one had formed a plausible conjecture for
$g(k)$ when $k>2$. Then, in the early 1990's, Conrey and
Ghosh~\cite{ConGho2} conjectured that $g(3) = 42$. Later, Conrey
and Gonek~\cite{ConGon} conjectured that $g(4) = 24024$. The
method employed by the last two authors reproduced the previous
values of $g(k)$ as well, but it did not produce a value for
$g(k)$ when $k>4$.

It was recently suggested by Keating and Snaith~\cite{KS} that the
characteristic polynomial of a large random unitary matrix can be
used to model the value distribution of the Riemann zeta function
near a large height $T$. Their idea was that because the zeta
function is analytic away from the point $s=1$, it can be
approximated at $s=\frac12 + \i t$ by polynomials whose zeros are
the same as the zeros of $\zeta(s)$ close to $t$. These zeros
(suitably renormalized) are believed to be distributed like the
eigenangles of unitary matrices chosen with Haar measure, so they
used the characteristic polynomial
\begin{equation}\label{Char poly}
Z_{N}(U,\theta) = \prod_{n=1}^{N} ( 1 - e^{\i(\theta_{n} - \theta)} )\,,
\end{equation}
where the $\theta_{n}$ are the eigenangles of a random $N \times
N$ unitary matrix $U$, to model $\zeta(s)$. For scaling reasons
they used matrices of size $N = \log T$ to model $\zeta(\frac12 + \i t)$
when $t$ is near $T$.
They then calculated the moments of
$|Z_{N}(U,\theta)|$ and found that
\begin{equation} \label{eq:mmts of char poly}
\E_N\left[ |Z_{N}(U,\theta)|^{2k} \right] \sim
\frac{G^{2}(k+1)}{G(2k+1)} N^{k^{2}}\,,
\end{equation}
where $\E_N$ denotes expectation with respect to Haar measure, and
$G(z)$ is Barnes' $G$-function. When $k = 1, 2, 3, 4$ they
observed that
\begin{equation*}
\frac{G^{2}(k+1)}{G(2k+1)} = \frac{g(k)}{\Gamma (k^2 +1)}\,,
\end{equation*}
where $g(k)$ is the same as in the results of Hardy--Littlewood
and Ingham, and in the conjectures of Conrey--Ghosh and
Conrey--Gonek given above.  They then conjectured that this holds
in general. That is, they asserted

\begin{conj}[Keating and Snaith]\label{conj:KS}
For $k$ fixed  with $\Re\,k>-1/2$,
\begin{equation*}
\frac{1}{T} \int_T^{2T} \left|\zeta(\tfrac{1}{2}+\i t)\right|^{2k}
\d t \sim a(k) \frac{G^2(k+1)}{G(2k+1)} (\log T)^{k^2} ,
\end{equation*}
as $T\to\infty$,
where $a(k)$ is given by \eqref{eq:a(k)}
and $G$ is the Barnes $G$--function.
\end{conj}

The characteristic polynomial approach has been successful in
providing insight into other important and previously intractable
problems in number theory as well (see, for example, \cite{MS} for
a survey of recent results). However, the model has the drawback
that it contains no arithmetical information---the prime numbers
never appear. Indeed, they must be inserted in an \emph{ad hoc}
manner. This is reflected, for example, by the absence of the
arithmetical factor $a(k)$ in equation \eqref{eq:mmts of char
poly}. Fortunately, in the moment problem it was only the factor
$g(k)$, and not $a(k)$, that proved elusive. A realistic model for
the zeta function (and other $L$-functions) clearly should include
the primes.

In this paper we present a new model for the zeta function that
overcomes this difficulty in a natural way. Our starting point is
an explicit formula connecting the zeros and the primes from
which we deduce a representation of the zeta function as a partial
Euler product times a partial Hadamard product.  Making certain
assumptions about how these products behave, we then reproduce
Conjecture \ref{conj:KS}.
Our model is based on the following representation of the zeta
function.

\begin{thm}\label{thm:zeta as product}
     Let $s=\sigma+\i t$ with $\sigma \geq 0$ and
     $ |t|\geq 2$, let $X \geq 2$ be a real parameter, and let $K$ be any fixed
     positive integer.
     Let $u(x)$ be a nonnegative $C^{\infty}$
     function of mass 1, supported on $[e^{1-1/X}, e]$, and set
\begin{equation}\label{U definition}
U(z) = \int_0^\infty u(x) E_1(z\log x)\;\d x \,,
\end{equation}
where $E_{1}(z)$ is the exponential integral
$\int_{z}^{\infty} e^{-w}/w\;\d w $.
Then
\begin{equation}\label{zeta formula}
    \zeta(s) = P_{X}(s) Z_{X}(s) \left(1+O\left(\frac{X^{K+2}}{(|s|
    \log X)^{K}}\right)+O(X^{-\sigma}\log X)\right)\;,
\end{equation}
where
\begin{equation} \label{eq:defn_P}
     P_{X}(s) = \exp\left(\sum_{n\leq X}\frac{\Lambda(n)}
     {n^{s}\log n} \right) \,,
\end{equation}
$\Lambda(n)$ is von Mangoldt's function, and
\begin{equation} \label{eq:defn Z_X}
   Z_{X}(s) = \exp\left(-\sum_{\rho_n}U\big((s-\rho_n)\log
    X\big)\right)\,.
\end{equation}
The constants implied by the $O$ terms depend only on $u$ and
$K$.
\end{thm}

We remark that Theorem~\ref{thm:zeta as product} is
unconditional---it does not depend on the assumption of any unproved
hypothesis. Moreover, it can easily be modified to accommodate
weight functions $u$ supported on the larger interval $[1, \,e]$.
Finally, as will be apparent from the proof, the second error term
can be deleted if we replace $P_{X}(s)$ by
\begin{equation*}
\widetilde P_{X}(s) =
   \exp\left(\sum_{n\leq X}\frac{\Lambda(n)}{n^{s}\log n}
   v(e^{\log n / \log X})\right)\,,
\end{equation*}
where $v(t) = \int_{t}^{\infty} u(x)\d x $.

To clarify \eqref{zeta formula}, we temporarily assume the Riemann
Hypothesis (RH) and take $s =\frac12+\i t$. We shall denote the
nontrivial zeros of $\zeta(s)$ by $\rho_n= \frac12 +\i\g_n $,
ordered by their height above the real axis, with $\g_{-n} =
-\g_n$. Since the support of $u$ is concentrated near $e$, $U(z)$
is roughly $E_{1}(z)$, which is asymptotic to $ -\gamma- \log z$
as $z \to 0$. Here $\gamma=0.5772\ldots$ is Euler's constant.
Thus, for those ordinates $\gamma_{n}$ close enough to $t$, we see
that
\begin{equation*}
\exp\,\big(-U(\i(t-\gamma_{n})\log X) \big)\approx
\i\,(t-\gamma_{n})\,e^{\gamma}\log X \,.
\end{equation*}
We expect the ordinates farther away not to contribute
substantially to the exponential defining $Z_{X}(s)$. Now,
$P_{X}(s) \approx \prod_{p \leq X}( 1 - p^{-s} )^{-1} $, hence our
formula looks roughly like
\begin{equation}\label{stripped down zeta formula}
\zeta(\tfrac12 +\i t) \approx \prod_{p \leq X}( 1 - p^{-\frac12-\i
t} )^{-1} \,
\prod_{\substack{\gamma_{n} \\
|t-\gamma_{n}|<1/\log X}} \bigg( \i(t-\gamma_{n})e^{\gamma}\log
X\bigg) \,.
\end{equation}

This formula is a ``hybrid'' consisting of a truncated Euler
product and (essentially) a truncated Hadamard product, with the
parameter $X$ mediating between them. Near height $T$ we are
approximating part of the zeta function by a polynomial of degree
about $ \log T/\log X$. The rest of the zeta function, which comes
from the zeros we have neglected, is approximated by the finite
Euler product. Formally, when we take $X$ large, we reduce the
number of zeros used to approximate zeta, but make up for it with
more primes; and when we take X small, we approach the previous
model~\eqref{Char poly}. Note however, that in order for the error
terms in \eqref{zeta formula} to be smaller than the main term, it
is necessary to work in an intermediate regime, where both the
zeros and the primes contribute.

To see how to use the model, and as a test case, we heuristically
calculate $I_{k}(T)$. The new model is more elaborate than the
original one, so more work is required. Nevertheless,
the idea is straightforward. The $2k$th moment of
$|\zeta(\frac12+\i t)|$ is asymptotic to the  $2k$th moment
of $|P_{X}(\frac12+\i t)\,Z_{X}(\frac12+\i t)|$.
We argue that when $X$ is not too large relative to $T$, the
$2k$th  moment of this product splits as the product of the moments.
We call this the ``Splitting Conjecture''.

\begin{conj}\label{thm:splitting}
{\bf (Splitting Conjecture.)} Let $X$ and $T\to \infty$ with $X
=O((\log T)^{2-\epsilon})$. Then for $k> -1/2$ we have
\begin{equation*}
\frac{1}{T} \int_T^{2T} \big|\zeta(\tfrac{1}{2}+\i t)\big|^{2k}
\d t \sim \bigg( \frac{1}{T}\int_T^{2T}
\big| P_{X}(\tfrac{1}{2}+\i t)\big|^{2k}\;
\d t  \bigg)\times
\bigg( \frac{1}{T}\int_T^{2T} \big|Z_{X}(\tfrac{1}{2}+\i t)\big|^{2k}
\; \d t \bigg)\,.
\end{equation*}
\end{conj}

In Section~\ref{sect:mean of P} we calculate the moments of $P$
rigorously and establish the following theorem.

\begin{thm}\label{thm:P moments}
Let  $1/2 \leq c <1$, $\epsilon>0$,
and let $k$ be any real number.
Suppose that $X$ and $T \to \infty$ and $X = O\left((\log
T)^{1/(1-c+\epsilon)}\right)$. Then we have
\begin{equation*}
\frac{1}{T}\int_T^{2T} |P_{X}(\sigma+\i t)|^{2k} \;\d t =
a(k,\sigma) F_{X}(k,\sigma) \left(1+ O_k\left(\frac{1}{\log X}
\right)\right)
\end{equation*}
uniformly for $c \leq \sigma\leq 1$,
where
\begin{equation}\label{eq:a(k,sigma)}
    a(k,\sigma) = \prod_{p }
    \left\{\left(1-\frac{1}{p^{2\sigma}}\right)^{k^2}
    \sum_{m=0}^{\infty}\frac{d_k(p^m)^{2}}{p^{2m\sigma}} \right\}
\end{equation}
and
\begin{equation*}\label{eq:F_{X}(k,sigma)}
F_{X}(k,\sigma ) =
\begin{cases}
\zeta(2\sigma)^{k^2} e^{-k^2 E_1((2\sigma-1)\log X)} & \text{ if }
\sigma>1/2 \;, \\
(e^{\gamma}  \log X)^{k^2} & \text{ if } \sigma=1/2 \;.
\end{cases}
\end{equation*}
Here $E_1$ is the exponential integral, and $\gamma=0.5772\ldots$
is Euler's constant.
\end{thm}

Note that $a(k, \frac12)$ is the same as  $a(k)$ in \eqref{eq:a(k)}.

In Section~\ref{sect:mean of Z} we  conjecture an asymptotic
estimate for $\int_{T}^{2T}|Z_{X}(\frac12+\i t)|^{2k}\d t$ using
random matrix theory. We introduce random matrix theory in the
following way. The statistical distribution of the ordinates
$\gamma_{n}$ is conjectured to coincide with that of the
eigenangles $\theta_{n}$ of $N \times N$ random unitary matrices
chosen with Haar measure for some $N$ (see for example~\cite{M},
\cite{O} and \cite{KeaSna3}). The choice of $N$ requires
consideration. The numbers $\gamma_{n}$ are spaced $2\pi /\log T$
apart on average, whereas the average spacing of the $\theta_{n}$
is $2 \pi/N$, and so we take $N$ to be the greatest integer less
than or equal to $\log T$. We therefore conjecture that the $2k$th
moment of $|Z_X(\frac12+\i t)|$, when averaged over $t$ around
$T$, is asymptotically the same as $|Z_X(\frac12+\i t)|^{2k}$ when
the $\gamma_n$ are replaced by $\t_n$ and averaged over all
unitary matrices with $N$ as specified above. We perform this
random matrix calculation in section~\ref{sect:mean of Z}, and so
obtain the following conjecture:

\begin{conj}\label{conj:RMT}
Suppose $X$, $T\to\infty$ with $X =O((\log T)^{2-\epsilon})$. Then
for any fixed $k> -1/2$, we have
\begin{equation*}
\frac{1}{T}\int_T^{2T} |Z_{X}(\tfrac12+\i t)|^{2k} \;\d t \sim
\frac{G^2(k+1)}{G(2k+1)}  \left(\frac{\log T}{ e^{\g }
\log X}\right)^{k^2}\;.
\end{equation*}
\end{conj}

We actually expect conjecture \ref{conj:RMT} to hold for a much
larger range of $X$, but the correct bound on the size of $X$ with
respect to $T$ is unclear.

We note that this asymptotic formula coincides with that in
(\ref{eq:mmts of char poly}) when there $N$ is taken to be on the
order of $\log T/e^{\gamma}\log X$.  This is consistent with the
fact that the polynomial in (\ref{stripped down zeta formula}) is
of about this degree. Alternatively, the mean density of
eigenvalues is $N$ divided by $2\pi$, and this is comparable to
the mean density of the ordinates of the zeros when multiplied by
$e^{\gamma}\log X$,  as they are in (\ref{stripped down zeta
formula}).

Combining the result of Theorem \ref{thm:P moments} with the
formula in Conjecture \ref{conj:RMT} and using the Splitting
Conjecture, we recover precisely the conjecture put forward by
Keating and Snaith.  Note that, as must be the case, all
$X$-dependent terms cancel out.

In Section~\ref{sect:splitting}, we prove

\begin{thm}\label{thm:first case of splitting}
     Let $\epsilon>0$ and let $X$ and $T \to \infty$
     with $X =O((\log T)^{2-\epsilon})$. Then for $k=1$ and $k=2$
     we have
     \begin{equation*}
     \frac{1}{T}\int_T^{2T}
     |\zeta(\tfrac12 + \i t ) P_{X}(\tfrac12 + \i t)^{-1} |^{2k } \;\d t
     \sim \frac{G^2(k+1)}{G(2k+1)}\,
     \left(\frac{\log T}{ e^{\gamma} \log X}\right)^{k^2} \,.
     \end{equation*}
\end{thm}

Since $ \zeta(\tfrac12 + \i t ) P_{X}(\tfrac12 + \i t)^{-1}
=Z_X(\tfrac12 + \i t) \left( 1+ o(1) \right)$ for $t \in[T, 2T]$,
it follows from this that Conjecture~\ref{conj:RMT} holds when $k=1$
and $k=2$. Moreover, combining Theorem~\ref{thm:first case of splitting}
with our estimate for
\begin{equation*}
\frac{1}{T}\int_T^{2T} |P_{X}(\tfrac12 + \i t)|^{2k } \;\d t
\end{equation*}
from Theorem~\ref{thm:P moments}, we also see that
Conjecture~\ref{thm:splitting} holds for $k=1$ and $k=2$.
Thus, we obtain the
\begin{cor*}
Conjectures \ref{thm:splitting} and \ref{conj:RMT} are true for
$k=1$ and $k=2$.
\end{cor*}

Clearly our model can be adapted straightforwardly to other
$L$-functions (see~\cite{KS1}). It can also be used to reproduce
other moment results and conjectures, such as those given by
Gonek~\cite{G} and by Hughes, Keating and O'Connell~\cite{HKO}
concerning derivatives of the Riemann zeta function at the zeros
of the zeta function. We also expect it to provide further insight
into the connection between prime numbers and the zeros of the
zeta function. It would be particularly interesting to determine
whether the model can be extended to capture lower order terms in
the asymptotic expansions of the moments of $\zeta(1/2+\i t)$ and
other $L$-functions, c.f.~\cite{CFKRS}.


\section{The Proof of Theorem 1}

We begin the proof by stating a
smoothed form of the explicit formula due to
Bombieri and Hejhal~\cite{BH}.

\begin{lem}\label{lem:BomHej}
Let $u(x)$ be a real, nonnegative, $C^{\infty}$ function with
compact support in $[1,e]$, and let $u$ be normalized so that if
\begin{equation*}
     v(t) = \int_t^{\infty} u(x) \;\d x \,
\end{equation*}
then $v(0)=1$. Let
\begin{equation*}
     \widetilde u(z) = \int_0^{\infty} u(x) x^{z-1} \d x
\end{equation*}
be the Mellin transform of $u$. Then for $s$ not a zero or pole of
the zeta function, we have
\begin{equation}\label{eq:explicit formula}
\begin{aligned}
     -\frac{\zeta'}{\zeta}(s) = &\sum_{n=2}^{\infty}
     \frac{\Lambda(n)}{n^{s}} v(e^{\log n/\log X})
     - \sum_{\rho} \frac{\widetilde u(1 - (s-\rho)\log X)}{s-\rho}
      + \frac{\widetilde u(1 - (s-1)\log X)}{s-1} \\
      -&\sum_{m=1}^{\infty} \frac{\widetilde u(1 - (s+2m)\log X)}{s+2m}
      \,,
\end{aligned}
\end{equation}
where the sum over $\rho$ runs over all the nontrivial zeros of the
zeta function.
\end{lem}

This lemma is proved in a familiar way, beginning with the
integral
\begin{equation*}
\frac{1}{2\pi\i}\int_{(c)} \frac{\zeta'}{\zeta}(z+s)
\widetilde{u}(1+z \log X)\frac{\d z}{z} \,,
\end{equation*}
where the integral is over the vertical line $\Re\, z = c =
\max\{2,2-\Re\ s\}$.

The support condition on $u$ implies that $v(e^{\log n/\log X})=0$
when $n>X$, so the sum over $n$ is finite. Furthermore, if $|\Im\,
z| >2$, say, then integrating $\widetilde{u}$ by parts $K$ times,
we see that
\begin{equation}\label{eq:estimate for Mellin transform}
    \begin{aligned}
     |\widetilde u(z)| \leq & \max_x | u^{(K)}(x)|  \,
     \bigg|\frac{\Gamma(z)}{\Gamma(z+K)} \bigg| \,
     (e^{\Re \, z + K} +1 )\\
  \leq & \max_x | u^{(K)}(x)| \,
      \frac{e^{\max\{\Re \, z + K,0\}} }{(1+|z|)^{K}}
\end{aligned}
      \end{equation}
for any positive integer $K$. Thus, the sums over $\rho$ and $m$
on the right-hand side of \eqref{eq:explicit formula} converge
absolutely so long as $s\neq \rho$ and $s\neq -2m$. This, in fact,
is the reason we require smoothing.

Next we integrate \eqref{eq:explicit formula} along the horizontal line
from  $s_{0} = \sigma_{0} + \i t_{0}$ to $+\infty$, where
$\sigma_{0} \geq 0$ and $|t_{0}| \geq 2$. If the line
does not pass through a zero, then on the left-hand side we obtain
$-\log\zeta(s_{0})$. We choose the branch of the logarithm here so that
$\lim_{\sigma \to \infty} \log \zeta(s) = 0$. If
the line of integration does pass through a zero, we define
$\log \zeta(\sigma + \i t) =  \lim_{\epsilon \to 0^{+}}
\frac{1}{2} \left(\log\zeta(\sigma + \i (t+ \epsilon)
+\log\zeta(\sigma + \i (t- \epsilon)\right)$\,.
Recalling the definition of $U(z)$ in \eqref{U definition},
we see that
\begin{equation}\label{Big U in terms of little u}
\begin{aligned}
\int_{s_{0}}^{\infty} \frac{\widetilde u(1-(s-z)\log X)}{s-z} \;\d s
&= \int_0^\infty u(x) E_1((s_{0}-z)\log X \log x) \;\d x\\
&= U((s_{0}-z)\log X) \,,
\end{aligned}
\end{equation}
provided that $s_{0} - z$ is not real and negative (so as to avoid
the branch cut of $E_1$).
If it is, we use the convention that
$U((s_{0}-z)\log X)=  \lim_{\epsilon \to 0^{+}}
\frac12 \big(U((s_{0}-z)\log X +\i \epsilon)  + U((s_{0}-z)\log X
-\i\epsilon)\big)$.
Note that
the logarithms in \eqref{Big U in terms of little u} are both positive
since the support of $u$ is in $[1,e]$ and $X \geq 2$.
It  therefore follows from \eqref{eq:explicit formula} that
\begin{equation}\label{eq:log zeta exact}
    \begin{aligned}
\log\zeta(s_{0}) = &\sum_{n=2}^\infty \frac{\Lambda(n)}{n^{s_{0}}\log n}
v(e^{\log n/ \log X}) - \sum_{\rho} U((s_{0}-\rho)\log X)\\
&+ U((s_{0}-1)\log X) -\sum_{m=1}^\infty U((s_{0}+2m)\log X) \,.
\end{aligned}
\end{equation}
The interchange of summation and integration in the sums
is justified by absolute convergence.
This representation holds for all points in $\Re\,s \geq 0$
not equal to the pole or one of the zeros of the zeta function.

We next suppose that the support of $u$  is contained in
$[e^{1-1/X},e]$ with the same $X$ as in \eqref{eq:log zeta exact}.
It is easy to see that there is a smooth nonnegative function $f$
with support in $[0,1]$ and total mass one such that $u(x) = X f(X
\log(x/e) + 1)/x$. Since $\max_x |f^{(K)}(x)|$ is bounded and
independent of $X$, we see that $\max_{x} |u^{(K)}(x)| \ll_{K}
X^{K+1}$. It therefore follows from \eqref{eq:estimate for Mellin
transform} that
\begin{equation*}
    \widetilde u(s) \ll_{K}
 \frac{e^{\max\{\sigma, 0\}} X^{K+1}}{(1+|s|)^{K}}\,.
\end{equation*}
{}From this and \eqref{Big U in terms of little u}, and
since $|t_{0}| \geq 2$, we find that
if $r$ is real, then
\begin{equation*}
\begin{aligned}
U((s_{0}-r)\log X) = &\int_{s_{0}}^{\infty}
\frac{\widetilde u(1-(s-r)\log X)}{s-r} \;\d s \\
&\ll_{K} \frac{X^{K+1}}{(\log X)^{K}}
\int_{\sigma_{0}}^\infty
\frac{X^{\max\{r-\sigma, \,0\}}}{|(\sigma-r) + \i t_{0}|^{K+1} }\;\d\sigma \\
&\ll_{K} \frac{X^{K+1+\max\{r-\sigma_{0}, \,0 \} }}{(\log X )^{K}}
\int_{\sigma_{0}}^\infty
\frac{1}{|(\sigma-r) + \i t_{0}|^{K+1} }\;\d\sigma \\
&\ll_{K} \frac{X^{K+1+\max\{r-\sigma_{0}, \,0 \}}}
{(|s_{0}-r|\log X)^{K}} \,.
\end{aligned}
\end{equation*}
In particular, for any fixed positive integer $K$ we have that
\begin{equation*}
U((s_{0}-1)\log X)
\ll_{K} \frac{X^{K+1+\max\{1-\sigma_{0}, \,0 \}}}
{(|s_{0}|\log X)^{K}} \;,
\end{equation*}
and, since $\sigma_{0} \geq 0$, that
\begin{equation*}
\begin{aligned}
\sum_{m=1}^\infty U((s_{0}+2m)\log X)
&\ll_{K} \frac{X^{K+1}}{(\log X)^{K}}
\sum_{m=1}^\infty \frac{1}{|s_{0}+2m|^{K+1}} \\
&\ll_{K} \frac{X^{K+1}}{(|s_{0}| \log X)^{K}}\;.
\end{aligned}
\end{equation*}

Inserting these estimates
into \eqref{eq:log zeta exact} and replacing $s_{0}$ by $s$,
we find that
\begin{equation*}
     \log\zeta(s) = \sum_{n=2}^\infty
     \frac{\Lambda(n)}{n^{s}\log n} v(e^{\log n/ \log X}) -\sum_{\rho}
     U((s-\rho)\log X) + O\left(\frac{X^{K+2}}{(|s|\log X)^K}\right)
\end{equation*}
for  $\sigma \geq 0$, $|t| \geq 2$,
and $K$ any fixed positive integer.
Exponentiating  both sides, we obtain
\begin{equation}\label{Temporary zeta product}
    \zeta(s) = \widetilde P_{X}(s) Z_{X}(s)
\left(1+ O\left(\frac{X^{K+2}}{(|s|\log X)^{K}}\right)\right)\;,
\end{equation}
where
\begin{equation*}
     \widetilde P_{X}(s) =
    \exp\left(\sum_{n\leq X}\frac{\Lambda(n)}{n^{s}\log n}
    v(e^{\log n / \log X})\right)
\end{equation*}
and
\begin{equation*}
  Z_{X}(s)= \exp\left(-\sum_{\rho} U((s-\rho) \log X)
\right)\,.
\end{equation*}

We now wish to show that replacing $\widetilde P_{X}(s)$ by
\begin{equation*}
  P_{X}(s) = \exp \left(\sum_{n\leq X}
\frac{\Lambda(n)}{n^{s}\log n} \right)
\end{equation*}
only introduces  a small error term into
\eqref{Temporary zeta product}.
To see this, note that $v((e^{\log n/\log X})) = 1$ for
$n\leq X^{1-1/X}$ because the support of $u(x)$ is
in $[e^{1-1/X},\,e]$. Therefore,
\begin{equation*}
\begin{aligned}
\frac{\widetilde P_{X}(s)}{P_{X}(s)} &= \exp\left(\sum_{X^{1-1/X}
\leq n \leq X}\frac{\Lambda(n)}{n^{s}\log n}
\left(v(e^{\log X/\log n}) - 1\right) \right) \\
 &\ll \exp\left(\sum_{X^{1-1/X} \leq n\leq X}\frac{1}{n^{\sigma}}
\right)\\
 &\ll \exp\left(X^{-\sigma}\log X\right) \;.
\end{aligned}
\end{equation*}

This completes the proof of Theorem~\ref{thm:zeta as product}
provided that $s$ is not a nontrivial zero of the zeta function.
To remove this restriction, we recall the formula
\begin{equation*}\label{Exp integral}
E_1(z) = - \log z  -\gamma  - \sum_{m=1}^{\infty}
\frac{(-1)^{m}z^{m}}{m!\,m} \,,
\end{equation*}
where $|\arg z| < \pi$, $\log z$ denotes the principal
branch of the logarithm, and $\gamma$ is Euler's constant.
{}From this and \eqref{U definition} we observe that we may interpret
$\exp(-U(z))$ to be asymptotic to $C z$ for some constant $C$
as $z\to 0$. Thus, both sides of \eqref{zeta formula}
vanish at the zeros.


\section{The Proof of Theorem~\ref{thm:P moments} } \label{sect:mean of P}

We begin with several lemmas.

\begin{lem}\label{lem:restricting_P_to primes and squares}
Let  $X \geq 2$ and set
\begin{equation*}\label{eq:P_ast}
P_X^\ast(s) = \prod_{p\leq X} \left(1-\frac{1}{p^s}\right)^{-1}
\prod_{\sqrt{X} < p \leq X} \left(1+   \frac{1}{2 p^{2s}}\right)^{-1}\,.
\end{equation*}
Then for $k$ any real number  we have
\begin{equation*}
P_X(s)^k = P_X^\ast(s)^k \left(1+O_k \left( \frac{1}{\log
X}\right)\right)
\end{equation*}
uniformly for $\sigma\geq 1/2$.
\end{lem}

\textbf{Proof.} By \eqref{eq:defn_P} we have
\begin{equation*}
P_{X}(s)^k = \exp\left(k\sum_{n\leq X}\frac{\Lambda(n)}{n^s \log n
} \right) =\prod_{p\leq X} \exp\left(k\sum_{1\leq j\leq N_p}
\frac{1}{j\,p^{j s}} \right) \,,
\end{equation*}
where $N_p = \left[\log X/\log p \right]$, the
integer part of $\log X / \log p$.
Therefore
\begin{align*}
P_{X}(s)^{k} P_X^\ast(s)^{-k} = \exp\left( -k\sum_{p\leq X}
\sum_{j> N_p} \frac{1}{j p^{js} } - k \sum_{\sqrt{X} < p \leq X}
\sum_{j=1}^\infty \frac{( -\tfrac12)^{j} }{j  p^{2js}} \right) \,.
\end{align*}
The primes $\sqrt{X} < p \leq X$ have $N_p =1$, and  we note that
the $j=2$ term for these primes in the first  double sum exactly
cancels the $j=1$ term in the second. Hence the argument in the exponent
is
\begin{align*}
& \ll |k|  \left(\sum_{p\leq \sqrt{X}}
\frac{1}{p^{\sigma (N_p+1)} }  +    \sum_{\sqrt{X} < p \leq X}
\frac{1}{p^{3\sigma}} \right)  \,.
\end{align*}
Now $p^{N_{p} +1} > X$ since  $N_{p} +1> \log X/\log p$,  so, for $\sigma \geq 1/2$,
this is
 \begin{align*}
& \ll  |k|\,  \left( X^{-1/2} \sum_{ p \leq \sqrt{X} }  1
  +    \sum_{ \sqrt{X} < p \leq X}
 \frac{1}{ p^{3/2}} \right)    \\
 & \ll    |k| \,  \left(\frac{1}{\log X}
 +  \frac{1}{X^{1/4}\log X} \right)
 \ll  \frac{|k|}{ \log X }\,.
\end{align*}
It follows that
 \begin{align*}
P_{X}(s)^{k} P_X^\ast(s)^{-k}
 = 1+ O_k \left(\frac{1}{\log X}\right) \,,
\end{align*}
as required.

\begin{lem}\label{lem:mean square dirichlet series}
Let $k$ be a real number.
Let $1/2 \leq c <1$ be arbitrary but  fixed,  and  suppose that
$2 \leq X \ll ( \log T)^{ 1/(1-c+\epsilon)}$, where  $\epsilon >0$ is
also fixed.
 Then
\begin{equation*}
\frac{1}{T} \int_T^{2T} |P_X^\ast(\sigma+\i t)|^{2k}\;\d t =
a(k,\sigma) \prod_{p\leq X}
\left(1-\frac{1}{p^{2\sigma}}\right)^{-k^2}  \left( 1 +
O_k\left(\frac{X^{1/2-2\sigma}}{\log X}\right) \right)
\end{equation*}
uniformly for $c \leq \sigma\leq 1$,
where $a(k,\sigma)$ is given by \eqref{eq:a(k,sigma)}.
\end{lem}

\textbf{Proof.}
We  write
\begin{equation}\label{P star to k}
\sum_{n=1}^\infty
\frac{\alpha_k(n)}{n^s} = P_X^*(s)^k =
\prod_{p\leq X} \left(1-p^{-s}\right)^{-k} \prod_{\sqrt{X} < p\leq
X} \left(1 +\tfrac12 p^{-2s}\right)^{-k}     \,.
\end{equation}
Let  $\smooth(X)$ denote the set of $X$-smooth numbers, that is,
$\smooth(X) = \{ n : p \mid n \implies p \leq X\}$.  Then
$\alpha_k(n) = d_k(n)$, the $k$th divisor function,  if $n \in
\smooth(\sqrt{X})$; $\alpha_k(p) = d_k(p)$ for all
$p\leq X$; and $\alpha_k(n) = 0$ if $n \not\in \smooth(X)$.
It is also easy to see that
\begin{equation*}
 \left(1-p^{-s}\right)^{-k}  \left(1 +\tfrac12 p^{-2s}\right)^{-k}
= \exp\left( k \left( \frac{1}{p^{s}} +  \frac{1}{3p^{3s}} +
\frac{1+ \left(\frac{1}{2}\right)^1}{4p^{4s}} + \frac{1}{5p^{5s}} +
 \frac{1- \left(\frac{1}{2}\right)^2}{6 p^{6s}}
+ \cdots     \right)     \right)  \,.
\end{equation*}
Comparing this with
\begin{equation*}
 \left(1-p^{-s}\right)^{-k}  =
 \exp\left( k \left( \frac{1}{p^{s}} +  \frac{1}{2p^{2s}} + \frac{1}{3p^{3s}} +
\frac{1 }{4p^{4s}}  + \cdots     \right)     \right)  \,,
\end{equation*}
we find that  for $k  \geq 0$, $\sqrt{X} < p \leq X$, and $j=1,2, \ldots$,
\begin{equation*}
 0 \leq \alpha_k(p^j) \leq d_{3k/2}(p^j) \,,
\end{equation*}
while for $k<0$
\begin{equation*}
 | \alpha_k(p^j)|  \leq   \alpha_{|k|}(p^j)  \leq d_{3|k|/2} (p^j) \,.
\end{equation*}

We now truncate the sum in \eqref{P star to k} at $T^\theta$, where $\theta$ is a small
positive number to be chosen later, and obtain
\begin{equation*}
\sum_{\substack{ n \in \smooth(X) \\ n \leq T^\theta}}
\frac{ \alpha_{k}(n) }{n^s} + O \left( \sum_{\substack{ n \in \smooth(X) \\
n > T^\theta}} \frac{ | \alpha_{k}(n) |}{n^{\sigma}} \right) \, .
\end{equation*}
For   $\epsilon  > 0$ fixed and $\sigma \geq c$,  the sum in the $O$-term is
\begin{align*}
& \ll   \sum_{\substack{ n> T^\theta\\ n\in\smooth(X) } }
\left(\frac{n}{T^{\theta} }\right)^\epsilon   \frac{d_{3|k|/2}(n)}{n^{\sigma}}
\leq T^{-\epsilon  \theta } \sum_{n \in \smooth(X)}
\frac{ d_{3|k|/2}(n)}{n^{c -\epsilon }} \\
&= T^{-\epsilon  \theta } \prod_{p\leq X}
\left(1 - p^{\epsilon  - c}   \right)^{-3|k|/2}
=  T^{-\epsilon  \theta }
\exp\left(  O\left( |k|\sum_{p\leq X}
p^{\epsilon  - c} \right)\right)  \\
& \ll   T^{-\epsilon  \theta }
\exp\left( O\left(\frac{|k|\,X^{ 1 - c + \epsilon   }}
{ (1- c + \epsilon ) \log X}\right)\right).
\end{align*}
Now suppose that $2 \leq X \ll (\log T)^{1/(1- c + \epsilon)}$
with the same $\epsilon$. Then  this is
\begin{align*}
 \ll   T^{-\epsilon \theta  }
\exp \left(  O\left(  \frac{ |k| \, \log T}
{ \log \log T}   \right)  \right)
  \ll_k   T^{-\epsilon \theta/2 } \,.
\end{align*}
Thus, we find that
\begin{equation}\label{eq:tmp_truncate_smooth_sum_in_P}
P_X^\ast(s)^k = \sum_{\substack{n \in  \smooth(X) \\ n \leq
T^\theta}} \frac{\alpha_k(n)}{n^s} + O_{ k}\left(T^{-\epsilon \theta/2}
\right) \,.
\end{equation}

Next we calculate  $ \frac{1}{T}\int_{T}^{2T}\left| P_X^\ast(s)  \right|^{2k} \d t$.
By  Montgomery and Vaughan's mean value theorem for Dirichlet polynomials~\cite{MV},
we have
\begin{align*}
 \int_T^{2T} \Biggl|\sum_{\substack{n \leq T^{\theta} \\
n \in \smooth(X)}} \frac{\alpha_k(n)}{n^{\sigma+\i t}}\Biggr|^2
 \;\d t
 &= \left( T + O( T^{\theta}  ) \right)
\sum_{\substack{n \leq T^{\theta}
\\ n \in \smooth(X)}} \frac{\alpha_k (n)^2}{n^{2\sigma}} \,.
\end{align*}
Using the method above, we may extend the sum on the right to infinity
with an error again no larger than $ O_{ k}\left(T^{-\epsilon \theta/2} \right)$.
Thus, taking $\theta = 1/2$, say, we find that
\begin{equation}\label{eq: mean of D poly}
\frac{1}{T} \int_T^{2T} \Biggl|\sum_{\substack{n \leq T^{1/2}\\
n \in \smooth(X)}} \frac{\alpha_k(n)}{n^{\sigma+\i t}}\Biggr|^2
 \;\d t
=  \sum_{n \in \smooth(X)} \frac{\alpha_k(n)^2}{n^{2\sigma}}
\left(1+O_{k}(T^{-\epsilon/4})\right) \, .
\end{equation}
We next note that if $A_i =  \frac{1}{T} \int_{T}^{2T} |a_i(t)|^2 \d t $ \,,
$i=1, 2 $, and $A_1 \neq 0$, then
$$
\frac1T\int_{T}^{2T} |a_1(t)+a_2(t)|^2 \d t =
 A_1  \left( 1 +  O\left( (A_2/A_1)^{1/2} ) \right) \right)  + A_2 \,,
$$
 the $O$-term arising from the Cauchy-Schwarz inequality applied
to the ``cross term".
We use this with $a_1(t)$  the sum on the right-hand side of
\eqref{eq:tmp_truncate_smooth_sum_in_P} and
$a_2(t)$ the error term  (with $\theta =1/2$).  Since $\alpha_k(1) =1$,
 we see from \eqref{eq: mean of D poly} that  $A_1 \gg 1$.
It therefore follows that
\begin{equation}\label{P* mean}
\frac{1}{T} \int_T^{2T} \left|P_{X}^\ast(\sigma +\i t)\right|^{2k}
\;\d t =
\left(1+O_{ k}(T^{-\epsilon/4})\right)
\sum_{n \in \smooth(X)} \frac{\alpha_k(n)^2}{n^{2\sigma}} \, .
\end{equation}
Since $\alpha_k(n) = d_k(n)$ for $n\in\smooth(\sqrt{X})$, and
$\alpha_k(p) = d_k(p)$ for $\sqrt{X} < p \leq X$, we may write
 the sum as
\begin{align*}
\sum_{n \in \smooth(X)} \frac{\alpha_k(n)^2}{n^{2\sigma}}
&=  \prod_{p \leq X} \left(  \sum_{m=0}^{\infty}
\frac{\alpha_k(p^m)^2}{p^{2m\sigma}} \right)  \\
& =  \prod_{p \leq \sqrt{X} } \left(  \sum_{m=0}^{\infty}
\frac{d_{k}(p^m)^2}{p^{2m\sigma}} \right)
 \prod_{\sqrt{X} < p \leq X} \left( 1 +  \frac{d_{k}(p)^2}{p^{2\sigma}}
 + \sum_{m=2}^{\infty} \frac{\alpha_k(p^m)^2}{p^{2m\sigma}} \right)  \\
&= \prod_{p \leq  X } \left(  \sum_{m=0}^{\infty}
\frac{d_{k}(p^m)^2}{p^{2m\sigma}} \right)
 \prod_{\sqrt{X} < p \leq X} \left\{
\left(  1 +  \frac{d_{k}(p)^2}{p^{2\sigma}}
 + \sum_{m=2}^{\infty} \frac{\alpha_k(p^m)^2}{p^{2m\sigma}}  \right)
\bigg{/} \sum_{m=0}^{\infty} \frac{d_{k}(p^m)^2}{p^{2m\sigma}}
  \right\} \;.
 \end{align*}
Factoring $1+ d_{k}(p)^2/p^{2\sigma}$  (which is at least $1$) out of the
numerator and denominator of the last product, we see that the product
equals
\begin{align*}
 \prod_{\sqrt{X} < p \leq X}    \left(   1 +
 O_{ k} \left( \frac{1}{p^{4\sigma}  } \right) \right)
 = \exp\left(O_{ k}\left( \frac{X^{\frac12 - 2\sigma}}
 {\log X} \right) \right)
  = 1+ O_{ k}\left( \frac{X^{\frac12 - 2\sigma}}{\log X} \right) \,.
 \end{align*}
Hence,
\begin{align*}
\sum_{n \in \smooth(X)} \frac{\alpha_k(n)^2}{n^{2\sigma}}
&=
\prod_{p \leq X} \left(  \sum_{m=0}^{\infty}
\frac{d_k (p^m)^2}{p^{2m\sigma}} \right)
 \left(  1+ O_{ k}\left(  \frac{X^{\frac12 - 2\sigma}}{\log X}  \right) \right) \,.
\end{align*}
Writing  the product here as
\begin{equation*}
\prod_{p\leq X} \left( \left(1-\frac{1}{p^{2\sigma}}\right)^{k^2}
\sum_{m=0}^{\infty} \frac{d_{k}(p^m)^2}{p^{2m\sigma}}\right)
\prod_{p\leq X} \left(1-\frac{1}{p^{2\sigma}}\right)^{-k^2} \, ,
\end{equation*}
we note that the first of the two factors may be extended over all the primes,
because
\begin{align*}
\prod_{p>X} \left(\left(1-\frac{1}{p^{2\sigma}}\right)^{k^2}
\sum_{m=0}^{\infty}\frac{d_k(p^m)^{2}}{p^{2m\sigma}}\right)
=& \prod_{p>X} \left(1+ O_{ k} \left(\frac{1}{p^{4\sigma}}\right) \right)  \\
= & 1+ O_{ k} \left(\frac{X^{1-4\sigma }}{\log X}\right) \,.
\end{align*}
Thus, by the definition of $a(k, \sigma)$  in \eqref{eq:a(k,sigma)}, we find that
\begin{align*}
\sum_{n \in \smooth(X)} \frac{\alpha_k(n)^2}{n^{2\sigma}}
 =
 a(k, \sigma) \prod_{p\leq X} \left(1-\frac{1}{p^{2\sigma}}\right)^{-k^2}
  \left(  1+ O_{  k} \left( \frac{ X^{\frac12 - 2\sigma}  }{\log X} \right) \right)
 \,.
\end{align*}
The lemma  follows from this and \eqref{P* mean}.

\begin{lem}\label{lem:the main term giving F}
If $k$ is a real number, then
\begin{equation*}
\prod_{p\leq X} \left( 1-\frac{1}{p^{2\sigma}} \right)^{-k^2} =
F_{X}(k,\sigma) \left(1+O_k(\frac{1}{\log X})\right)
\end{equation*}
uniformly for $\sigma \geq 1/2$, where
\begin{equation*}
F_{X}(k,\sigma) =
\begin{cases}
\zeta(2\sigma)^{k^2} e^{-k^2 E_1((2\sigma-1)\log X)} & \text{ if }
\sigma>1/2 \,,\\
(e^{\gamma } \log X)^{k^2} & \text{ if } \sigma=1/2 \,,
\end{cases}
\end{equation*}
and $E_1$ is the exponential integral.
\end{lem}

\textbf{Proof.}  Mertens'  theorem asserts that
\begin{equation*}
     \prod_{p\leq X} \left( 1-\frac{1}{p} \right)^{-1} =
     e^{\gamma}\log X \left(1+O(\frac{1}{\log X})\right) \,.
\end{equation*}
Raising both sides to the $k^{2}$ power establishes the result
when $\sigma=1/2$. When $\sigma>1/2$, we see that
\begin{equation*}
     \prod_{p\leq X} \left( 1-\frac{1}{p^{2\sigma}} \right)^{-1} =
     \zeta(2\sigma) \exp\left(\sum_{p>X}
     \log\left(1-\frac{1}{p^{2\sigma}}\right)\right) \,.
\end{equation*}
By the prime number theorem in the form $\psi(x)=\sum_{n\leq
x}\Lambda(n) = x + O(x/(\log x)^{A})$, we find that
\begin{align*}
     \sum_{p>X}\log\left(1-\frac{1}{p^{2\sigma}}\right)
     &= -\sum_{p>X}\left( \frac{1}{p^{2\sigma}} +
     O\left(\frac{1}{p^{4\sigma}}\right) \right) \\
     &= -\int_{X}^{\infty} \left(\frac{1}{u^{2\sigma}} +
     O\left(\frac{1}{u^{4\sigma}}\right)\right) \frac{\d u}{\log u}
     + O\left(\frac{1}{(\log X)^{A}}\right) \\
     &=-E_1((2\sigma-1)\log X) + O\left(\frac{1}{(\log X)^A}\right)\,.
\end{align*}
Hence,
\begin{equation*}
     \prod_{p\leq X} \left( 1-\frac{1}{p^{2\sigma}} \right)^{-k^{2}} =
     \zeta(2\sigma)^{k^{2}} \exp\left(-k^{2} E_1((2\sigma-1)\log X)\right)
     \left( 1+ O_{k}\left(\frac{1}{\log X}\right) \right)\,,
\end{equation*}
as asserted.

The proof of Theorem~\ref{thm:P moments} now follows immediately
from Lemmas~\ref{lem:restricting_P_to primes and squares},
\ref{lem:mean square dirichlet series} and
\ref{lem:the main term giving F}.


\section{Support for Conjecture~\ref{conj:RMT}}
\label{sect:mean of Z}

In this section we give heuristic arguments supporting
Conjecture~\ref{conj:RMT}, which we  restate as
\begin{equation*}
\frac{1}{T} \int_T^{2T} \left|Z_X(\tfrac12 + \i t)\right|^{2k}
\;\d t \sim \frac{G^2(k+1)}{G(2k+1)}\,
\left(\frac{\log T}{ e^{\g }\log X}\right)^{k^2}
\end{equation*}
as $T \rightarrow \infty$, where $Z_X(s)$ is given by
\eqref{eq:defn Z_X}.

We assume the Riemann Hypothesis. Since $\Re\, E_1(\i x) =
-\Ci(|x|)$ for $x\in\R$, where
\begin{equation*}
\Ci(z)= -\int_{z}^{\infty} \frac{\cos w}{w}\,\d w \,,
\end{equation*}
we find that
\begin{equation}\label{eq:defn int |Z_X|}
\frac{1}{T}\int_T^{2T} \left|Z_X(\tfrac12 + \i t)\right|^{2k} \;\d t =
\frac{1}{T}\int_T^{2T} \prod_{\gamma_n} \exp\left( 2k \int_1^e u(y)
\Ci(|t-\gamma_n|\log y\log X)\;\d y \right) \;\d t \,,
\end{equation}
where $u(y)$ is a smooth nonnegative function, supported on
$[e^{1-1/X},e]$ and of total mass $1$. Since the terms in the
exponent decay as $|\gamma_n - t |$ increases, this product is
effectively a local statistic. That is, the integrand depends only
on those zeros close to $t$. In recent years considerable evidence
has been amassed suggesting that the zeros of the Riemann zeta
function around height $T$ are distributed like the  eigenangles
of unitary matrices of size $\log T$ chosen with Haar measure
(see, for example, the survey article \cite{KeaSna3}). We
therefore model the right-hand side of \eqref{eq:defn int |Z_X|}
by replacing the ordinates $\gamma_n$ by the eigenangles of an $N
\times N$ unitary matrix and averaging over all such matrices with
Haar measure, where $N=[\log T]$. Thus, the right-hand side of
\eqref{eq:defn int |Z_X|} should be asymptotic to
\begin{equation*}
\E_N\left[\prod_{n=1}^N \exp\left( 2k \int_1^e u(y)
\Ci(|\theta_n|\log y\log X)\;\d y \right) \right] \,,
\end{equation*}
where the $\t_n$ are the eigenangles of the random matrix and
$\E_N\left[ \cdot \right]$ denotes the expectation with respect to Haar measure.
However, since the eigenangles of a unitary matrix are
naturally $2\pi$-periodic objects, it is convenient to periodicize
our function, which we do by defining
\begin{equation}\label{eq:symbol}
\phi(\t) =  \exp\left( 2k  \int_1^e u(y) \left( \sum_{j=-\infty}^\infty
\Ci(|\theta+2\pi j|\log y\log X)\right) \;\d y \right) \,.
\end{equation}
It will follow from our proof of Lemma~\ref{lem:b(0)} that the
 terms with $j \neq 0$, which make the random matrix
calculation much easier, only contribute $\ll_k 1/\log X$ to
$\phi(\t)$ when $-\pi < \t \leq \pi$. Hence they do not affect the
accuracy of the model.
Thus, we argue that
\begin{equation}\label{eq:RMT model}
\frac{1}{T}\int_T^{2T} \left|Z_X(\tfrac12 + \i t)\right|^{2k} \;\d
t \sim \E_N \left[ \prod_{n=1}^N \phi(\t_n) \right] \,.
\end{equation}
The remainder of this section is devoted to the proof of
\begin{thm}\label{thm:RMT model}
Let $\phi(\t)$ be defined as in \eqref{eq:symbol}, then for fixed
$k>-\frac12$ and $X\geq 2$, we have as $N\to\infty$,
\begin{equation*}
\E_N \left[ \prod_{n=1}^N \phi(\t_n) \right] \sim
\frac{(G(k+1))^2}{G(2k+1)} \left(\frac{N}{e^\gamma \log X }
\right)^{k^2} \left(1+O_k\left(\frac{1}{\log X}\right)\right) \,.
\end{equation*}
\end{thm}

\begin{rem}
The random matrix model of Keating and Snaith \cite{KS} for the
moments of the Riemann zeta function involved the characteristic
polynomial \eqref{Char poly}. Note that if we set
$M=\frac{N}{e^\gamma \log X}$, then by \eqref{eq:mmts of char poly}
we have
\begin{equation*}
\E_{M} \left[\left|Z_{M}(U,\t)\right|^{2k} \right] \sim
\frac{(G(k+1))^2}{G(2k+1)}\left(\frac{N}{e^\gamma \log X }
\right)^{k^2} \,,
\end{equation*}
which is the same answer we find in Theorem~\ref{thm:RMT model}.
This is easily explained by the fact that in our model the
eigenangles are multiplied by $e^\gamma \log X$ and so their mean
density is $M/2\pi$.  Given that for random matrices the mean
density is the only parameter in the asymptotics of local
eigenvalue statistics, it is natural  that the result should
be the same as for unitary matrices of dimension $M$, since their
eigenangles have precisely this mean density.
\end{rem}

\textbf{Proof.}   Heine's identity
\cite{szego} evaluates the expected value in \eqref{eq:RMT model}
as a Toeplitz determinant
\begin{equation}\label{eq:Heine}
\E_N \left[ \prod_{n=1}^N \phi(\t_n) \right] = \det \left[
\phi_{i-j} \right]_{1\leq i,j\leq N} \ ,
\end{equation}
where
\begin{equation*}
\phi_n = \frac{1}{2\pi} \int_{-\pi}^\pi \phi(\t) e^{-\i n \t}
\;\d \t
\end{equation*}
is  the $n$th Fourier coefficient  of $\phi(\t)$. The Toeplitz  symbol
$ \phi(\t)$ is singular since it is zero when $\t=0$. Thus,
 the asymptotic evaluation of this determinant requires
knowledge of the Fisher--Hartwig Conjecture in a form proved by
Basor \cite{B}.

We factor out the singularity in $\phi(\t)$ by writing
\begin{equation*}
\phi(\t) = b(\t) (2-2\cos \t)^k\,,
\end{equation*}
where
\begin{equation}\label{eq:b}
b(\t) = \exp\left(-k\log(2-2\cos \t) +  2k \int_1^e u(y) \left( \sum_{j=-\infty}^\infty
 \Ci(|\theta+2\pi j|\log y\log X) \right) \;\d y \right) \,.
\end{equation}
As we will see in the proof of Lemma~\ref{lem:b(0)} below,  the
logarithmic singularities in the exponent on the right cancel.
Thus $b(\t)$ never equals zero. The asymptotic behavior of the
Toeplitz determinant with these symbols has been determined by
Basor \cite{B}. She showed that if $k>-1/2$, then
\begin{equation}\label{eq:Fisher-Hartwig}
\det \left[ \phi_{i-j} \right]_{1\leq i,j\leq N} \sim E
\exp\left(\frac{N}{2\pi} \int_{-\pi}^\pi \log b(\t) \;\d \t
\right) N^{k^2}
\end{equation}
as $N\to\infty$, where the constant $E$ is given by
\begin{equation*}
E = \exp\left(    \sum_{n=1}^\infty  n \left(\frac{1}{2\pi}
\int_{-\pi}^\pi \log b(\t) e^{- \i n \t} \;\d \t\right)^2 \right)
b(0)^{-k} \frac{ G^2(k+1)}{G(2k+1)} \,.
\end{equation*}

To evaluate $E$ we need to know $b(0)$  and the Fourier
coefficients of $\log b(\t)$. These are given by the next two
lemmas.

\begin{lem}\label{lem:fourier coeff log b}
Let $b(\t)$ be given by \eqref{eq:b}. Then
\begin{equation*}
\frac{1}{2\pi} \int_{-\pi}^\pi \log b(\t) e^{-\i n\t}\;\d \t =
\begin{cases}
0 & \text{ if } n = 0 \,, \\
\frac{k}{n} v \left(e^{n/\log X}\right) & \text{ if } n\geq 1 \,,
\end{cases}
\end{equation*}
where
\begin{equation*}
v(t) = \int_t^\infty u(y)\;\d y\,.
\end{equation*}
\end{lem}

\begin{lem}\label{lem:b(0)}
Let $b(\t)$ be given by \eqref{eq:b} and let $u(x)$ have total mass
one with support in $[e^{1-1/X},e]$. Then
\begin{equation*}
b(0) =\exp \bigg(  2k \left( \log\log X + \gamma \right) \bigg)
\left( 1 + O_k \left( \frac{1}{\log X} \right)  \right) \,.
\end{equation*}
\end{lem}

Before proving the lemmas, we complete the proof of
Theorem~\ref{thm:RMT model}. Since $u$ is a nonnegative function
supported in $[e^{1-1/X},e]$ of total mass one, we see that
\begin{equation*}
v(t) =
\begin{cases}
1 & \text{ if } t\leq e^{1-1/X} \,,\\
0 & \text{ if } t \geq e \,,
\end{cases}
\end{equation*}
and $ 0 \leq v(t) \leq 1$ if $  t \in [e,  e^{1-1/X}]$. Thus
\begin{align*}
\sum_{n=1}^\infty \frac{1}{n} \left(v \left(\exp(\frac{n}{\log
X})\right)\right)^2
  = \sum_{n\leq (1-1/X)\log X } \frac{1}{n}
  + O\left( \sum_{(1-1/X) \log X<n \leq \log X}  \frac{1}{n} \right)\,.
\end{align*}
The first sum on the right equals $\log\log X + \gamma
+O\left( 1/\log X \right)$ and the second is $O\left(X^{-1}\right)$ .
Hence, we find that
\begin{equation*}
\sum_{n=1}^\infty \frac{1}{n} v \left(\exp(\frac{n}{\log
X})\right)^2 = \log\log X + \gamma +O\left(\frac{1}{\log X}\right) \,.
\end{equation*}
Using this and the value of $b(0)$ given by Lemma~\ref{lem:b(0)},
we obtain
\begin{align*}
E =\exp\left( -k^2 (\log \log X + \gamma) \right)
\frac{(G(k+1))^2}{G(2k+1)} \left(1+O_k\left(\frac{1}{\log
X}\right)\right) \,.
\end{align*}

The proof of Theorem~\ref{thm:RMT model} is completed by
combining this, the case $n=0$ of  Lemma~\ref{lem:fourier coeff
log b},   \eqref{eq:Heine}, and  \eqref{eq:Fisher-Hartwig}.

\textbf{Proof of Lemma~\ref{lem:fourier coeff log b}.} We wish to
evaluate
\begin{equation*}
\frac{1}{2\pi} \int_{-\pi}^\pi \log b(\t) e^{-\i n\t}\;\d \t \,,
\end{equation*}
where $b(\t)$ is given by \eqref{eq:b}. After some straightforward algebra
we see that this equals
\begin{equation}\label{eq:FC log b}
\frac{-k}{\pi} \int_0^\pi \log(2-2\cos \t) \cos n\t \;\d \t +
\frac{2k}{\pi}   \int_1^e  u(y) \left( \int_{0}^\infty\Ci(\theta \log y\log
X) \cos n \t  \;\d\t \right) \;\d y \,.
\end{equation}

When $n=0$ the first integral vanishes by symmetry,
and the second vanishes because
\begin{align*}
  \int_{0}^\infty   \Ci(\t) \;\d \t = 0 \;.
\end{align*}
This is a special case of  the  formula (see Gradshteyn and Ryzhik~\cite{GR}, p. 645)
\begin{align}\label{cosine int}
\int_{0}^\infty \Ci( A \theta  ) \cos n \t \;\d\t  =
\begin{cases}
-\frac{\pi}{2n} & \text{ if } A<n,\\
-\frac{\pi}{4n} & \text{ if } A=n,\\
0 & \text{ otherwise}
\end{cases}
\end{align}
for $A>0$, which we require below as well. Thus, both terms in
\eqref{eq:FC log b} vanish and Lemma~\ref{lem:fourier coeff log b} holds in this case.

When $n $ is a positive integer, the first term in \eqref{eq:FC log b}
equals
\begin{align}\label{eq:log int FC}
-\frac{ k}{\pi} \int_{0}^\pi \log (2-2\cos \t ) \cos n\t \;\d \t &=
-\frac{ k}{\pi} \int_{0}^\pi \left(\log 4 + 2  \log (\sin \frac{\t}{2} ) \right) \cos n\t \;\d \t  \nonumber \\
&=  -\frac{ 4k}{\pi}  \int_{0}^{\pi/2 }  \log (\sin \t)   \cos 2n \t \;\d \t  \\
&= \frac{ k}{n}    \nonumber
\end{align}
(see Gradshteyn and Ryzhik~\cite{GR}, p. 584).
The second term in \eqref{eq:FC log b}  is, by \eqref{cosine int},
\begin{align}
\frac{2k}{\pi} \int_1^e  u(y) \left( \int_{0}^\infty
\Ci(\theta \log y\log X) \cos n \t \;\d\t  \right )\;\d y
&= -\frac{k}{n} \int_{1}^{e^{n/\log X}} u(y) \;\d y \nonumber\\
&=  \frac{k}{n} \left(v(e^{n/\log X}) -1\right) \,.\label{eq:FC
cosine int}
\end{align}
Inserting \eqref{eq:FC cosine int} and \eqref{eq:log int FC} into
\eqref{eq:FC log b}, we find that for $n> 0$ an integer,
\begin{equation*}
\frac{1}{2\pi} \int_{-\pi}^\pi \log b(\t) e^{-\i n\t}\;\d \t =
\frac{k}{n} v(e^{n/\log X})  \,.
\end{equation*}
This completes the proof of Lemma~\ref{lem:fourier coeff log b}.

\textbf{Proof of Lemma~\ref{lem:b(0)}.} We  calculate
$b(0)$, where
\begin{equation}\label{b theta}
b(\t) = \exp\left(-k\log(2-2\cos \t) +  2k \int_1^e u(y)  \left( \sum_{j=-\infty}^\infty
\Ci(|\theta+2\pi j|\log y\log X)\right)  \;\d y \right) \,.
\end{equation}
Using the expansion
\begin{equation*}
\Ci(x) = \gamma + \log x + O(x^2)
\end{equation*}
for $x>0$, we find that the first term in the exponent and the $j=0$ term combined
contribute
\begin{multline*}
-k\log(2-2\cos \t) +  2k \int_1^e u(y) \Ci(|\theta|\log y\log X)\;\d y   \\
=  2k \left\{ - \log (|\t|)  + O(\t^2)  +   \int_1^e u(y)
\left( \log(|\theta|\log y\log X) + \gamma + O_X(\t^2) \right)
 \;\d y    \right\}     \\
=   2k \left\{ \gamma + \log\log X  +  \int_1^e u(y) \log\log y \;\d y +
O_X(\t^2)  \right\} \,,
\end{multline*}
since $u(x)$ has total mass one.
Moreover,  $u(x)$ is supported in $[e^{1-1/X}, e]$, so we have
\begin{equation*}
\int_1^e u(y) \log\log y \;\d y  \ll \frac{1}{X} \,.
\end{equation*}
Therefore we find that
\begin{align}\label{part of b(0)}
\lim_{\t \to 0}  \bigg\{ -k \log(2-2\cos \t) &+  2k \int_1^e u(y)
\Ci( | \theta | \log y \log X) \;\d y   \bigg\} \\
&=  2k \left( \log\log X + \gamma \right) +
O_k \left( \frac{1}{X} \right) \,.  \notag
\end{align}

Now consider the contribution of the terms with $j \neq 0$ in
\eqref{b theta}. An integration by parts shows that
\begin{equation*}
\Ci(x) = -\int_{x}^{\infty} \frac{\cos t}{t}\,\d t
= \frac{\sin x}{x} + O(\frac{1}{x^2})
\end{equation*}
for $x$ positive and $  \gg 1$. Thus, since $(1-1/X) \log X
\leq  \log y \; \log X \leq \log X$, $X >2$, and $\t\in(-\pi,\pi]$, we see that
\begin{align*}
 \sum_{\substack{ j=-\infty \\ j\neq 0}}^{\infty}
 \Ci \left( |\theta  + 2 \pi j | \log y \log X \right)
 &=  \frac{1}{\log y\log X} \sum_{\substack{ j=-\infty\\ j\neq 0}}^{\infty}
\frac{ \sin \left(|\theta+2\pi j | \log y\log X \right) }{|\theta+2\pi j | } \;
+ \; O\left(  \frac{1}{ (\log X)^2 } \right) \;.
\end{align*}
In a standard way (\emph{via} Abel partial summation), one can show that the series on the right is uniformly
convergent for $y \in [e^{1-1/X}, e]$,  except possibly  in the neighborhood
of a finite number of points, and boundedly convergent over the whole interval. Moreover,
 the series may be bounded  independently  of $\t \in (-\pi, \pi]$.
We may therefore  multiply by the continuous function $u(y)$ and integrate to find that
 \begin{equation*}
\int_1^e u(y) \left( \sum_{\substack{j=-\infty\\ j \neq 0}}^\infty
\Ci(|\theta+2\pi j|\log y\log X) \right) \;\d y
\ll \frac{1}{\log X} \int_1^e \frac{ u(y)}{\log y}   \;\d y
+  O\left( \frac{1}{(\log X)^2}  \right)
\ll \frac{1}{\log X} \,.
\end{equation*}
uniformly for  $\t \in (-\pi, \pi ]$.
Combining this and \eqref{part of b(0)}  with \eqref{b theta}, we obtain
\begin{align*}
 b(0) &= \exp \left(  2k \left( \log\log X + \gamma \right) +
O_k \left( \frac{1}{\log X} \right)  \right) \\
&= \exp \bigg(  2k \left( \log\log X + \gamma \right) \bigg)
\left( 1 + O_k \left( \frac{1}{\log X} \right)  \right)\,.
\end{align*}
This completes the proof of Lemma~\ref{lem:b(0)}.


\section{The Proof of Theorem 3}\label{sect:splitting}

First we prove Theorem~\ref{thm:first case of splitting} when
$k=1$. In this case $G^2(k+1)/G(2k+1)= G^2(2)/G(3)=1$, and by
Lemma~\ref{lem:restricting_P_to primes and squares} we may replace
$ P_X(\tfrac12+\i t)$ by $ P_X^\ast(\tfrac12+\i t)$. Thus, it
suffices to show that for $X \ll (\log T)^{2-\epsilon}$,
\begin{equation*}
\frac{1}{T} \int_T^{2T} \left|
 \zeta(\tfrac12+\i t) P_X^\ast(\tfrac12+\i t)^{-1} \right|^{2}
 \;\d t =
\frac{\log T}{ e^{\gamma} \log X} \left(1+ O\left(\frac{1}{\log X}\right)\right)\,.
\end{equation*}

As in the proof of   Lemma~\ref{lem:mean square dirichlet series}, we write
 $\smooth(X) = \{ n :  p \mid n \implies p\leq X\}$ and
\begin{equation*}
P_X^\ast(\tfrac12+\i t)^{-1} = \sum_{n\in \smooth(X)}
\frac{\alpha_{-1}(n)}{n^{1/2 + \i t}} \,,
\end{equation*}
where $\alpha_{-1}(n) = \mu(n)$, the M\"obius function, if $n\in
\smooth(\sqrt{X})$; $\alpha_{-1}(p) = \mu(p)$ for all $p\leq X$; and
$\alpha_{-1}(n) \ll d_{3/2}(n) \ll d(n)$ for all $n\in \smooth(X)$.
By \eqref{eq:tmp_truncate_smooth_sum_in_P},
if the $\epsilon$  above is sufficiently small, we  find that
\begin{equation}\label{eq:ranking splitting the sum 1}
P_X^\ast(\tfrac12 +\i t)^{-1} = \sum_{\substack{n \leq T^\theta \\ n \in
\smooth(X)}} \frac{\alpha_{-1}(n)}{n^{1/2+\i t}} +
O\left(T^{-\theta\epsilon/10} \right)
\end{equation}
(The exponent $1/10$ in place of $1/2$ is accounted for by the slight difference between
 the conditions $X \ll (\log T)^{2-\epsilon}$ and $X \ll (\log T)^{1/(1/2 + \epsilon)}$.)
Now for $m$ and $n$  coprime positive integers, we have the   formula
\begin{align*}
 \int_T^{2T} \left| \zeta(\tfrac12+\i t) \right|^2
\left( \frac{m}{n} \right)^{\i t} \;\d t
=  \frac{T}{ \sqrt{mn} } \left( \log\left( \frac{T}{2\pi mn }\right) + 2\gamma -1 \right)
+  O\left( mn T^{8/9} (\log T)^6 \right) \,.
\end{align*}
(For example, see Corollary~24.5 of \cite{IK}.)
Using this and the main term in \eqref{eq:ranking splitting the sum 1} with
$\theta=1/20$, we find that
\begin{multline}\label{int formula}
\frac{1}{T} \int_T^{2T} \left|\zeta(\tfrac12+\i t)\right|^2 \Bigg|
\sum_{\substack{n \leq T^{1/20} \\ n \in \smooth(X)}}
\frac{\alpha_{-1}(n)}{n^{1/2+\i t}} \Bigg|^2 \;\d t \\
= \sum_{\substack{m,n \leq T^{1/20} \\ m,n \in \smooth(X)}}
\frac{\alpha_{-1}(m)}{m} \frac{\alpha_{-1}(n)}{n} (m,n) \left\{\log
\left(\frac{T(m,n)^2}{2\pi m n}\right) + 2\gamma -1 +
O\left(\frac{m n}{(m,n)^2} T^{ -1/9}  (\log T)^6  \right)\right\}\,,
\end{multline}
where  $(m,n)$ denotes the greatest common divisor of $m$ and $n$.
The $O$-term contributes
\begin{equation*}
\ll  T^{-1/9}  (\log T)^6 \left(\sum_{n  \leq T^{1/20}}
d (n)  \right)^2
\ll T^{ -1/90}  (\log T)^8 \,.
\end{equation*}
Grouping together those $m$ and $n$ for which
$(m, n)=g$, replacing $m$ by $gm$ and $n$ by $gn$,
and then using the inequality $d(ab)\leq d(a)d(b)$, we find that
\begin{multline*}
\sum_{\substack{m,n \leq T^{\theta} \\ m,n \in \smooth(X)}}
\frac{\alpha_{-1}(m)}{m} \frac{\alpha_{-1}(n)}{n} (m,n) \left(\log
\left(\frac{(m, n)^2}{2\pi m n}\right) + 2\gamma -1\right) \\
\ll \sum_{g \in\smooth(X)} \frac{1}{g} \sum_{\substack{m, n
\in\smooth(X) \\ (m, n)=1}} \frac{d(gm) d(gn)\log mn}{mn}    \ll
\sum_{g \in\smooth(X)} \frac{d(g)^2}{g}\left(  \sum_{n \in
\smooth(X)} \frac{d(n) \log n}{n} \right)^2 \,.
\end{multline*}
If we write $f(\sigma)= \sum_{n \in\smooth(X)} d(n)n^{-\sigma} =
\prod_{p \leq X} \left(1- p^{-\sigma}\right)^{-2}$, then the sum
over $n$ is $-f^{'}(1)$, which, by logarithmic differentiation, is
$ 2 f(1) \sum_{p \leq X} \log p/ (p-1) \ll f(1) (\log X) \ll (\log
X)^3.$ We also have $\sum_{g \in\smooth(X)} d(g)^2 g^{-1} \ll
\prod_{p \leq X} \left(1- p^{-1}\right)^{-4} \ll  (\log X)^4$, and
so the expression above is $\ll (\log X)^{10}$.

Thus far then, we have
\begin{align}\label{int formula 2}
\frac{1}{T} \int_T^{2T} \left|\zeta(\tfrac12+\i t)\right|^2 &\Bigg|
\sum_{\substack{n \leq T^{1/20} \\ n \in \smooth(X)}}
\frac{\alpha_{-1}(n)}{n^{1/2+\i t}} \Bigg|^2 \;\d t   \\
& =  \log T \sum_{\substack{m,n \leq T^{1/20} \\ m,n \in \smooth(X)}}
\frac{\alpha_{-1}(m)}{m} \frac{\alpha_{-1}(n)}{n} (m,n)
 \quad + O\left(  (\log X)^{10} \right) \,. \notag
\end{align}
Since  $\sum_{g | n}  \phi(g) = n$, the remaining sum here  is
\begin{equation}\label{main term sum 1}
\sum_{ \substack{m,n \leq T^{1/20} \\  m, n \in \smooth(X)} }
 \frac{\alpha_{-1}(m)}{m} \frac{\alpha_{-1}(n)}{n}
  \left( \sum_{\substack{g|m\\ g|n}} \phi(g)  \right)
  = \sum_{ \substack{g \leq T^{1/20} \\  g \in \smooth(X)}} \frac{\phi(g)}{g^2}
 \left(    \sum_{ \substack{ n \leq  T^{1/20}g^{-1} \\    n \in \smooth(X)} }
\frac{\alpha_{-1}(gn)}{n} \right)^2\,.
\end{equation}
We wish to extend the  sums on the right  to all of $\smooth(X)$.
For this we use several  estimates.
First,
$$
\sum_{n \in \smooth(X)} \frac{|\alpha_{-1}(gn)|}{n}
\leq d(g) \sum_{n \in \smooth(X)} \frac{d(n)}{n}
=d(g) \prod_{p \leq X}\left( 1-\frac{1}{p} \right)^{-2}
\ll d(g)( \log X)^2 \,.
$$
Second,
\begin{align*}
\sum_{\substack{ n> T^{1/20}g^{-1} \\ n \in \smooth(X)} }\frac{|\alpha_{-1}(gn)|}{n}
&\leq d(g) \sum_{\substack{  n > T^{1/20}g^{-1}    \\n \in \smooth(X)}} \frac{d(n)}{n}
\leq d(g) \left( \frac{T^{1/20}}{g} \right)^{-1/4} \sum_{n \in \smooth(X)} \frac{d(n)}{n^{3/4}} \\
&\ll d(g)g^{ 1/4} T^{-1/80} \prod_{p\leq X}\left(  1- \frac{1}{p^{3/4}}\right)^{-2}
\ll d(g)g^{ 1/4} T^{-1/80} e^{10 X^{1/4}/\log X} \\
& \ll d(g)g^{ 1/4} T^{-1/100} \,,
\end{align*}
say. From these it follows that the square of the sum over $n$ in \eqref{main term sum 1} is
\begin{align}\label{sum squared}
 \left(    \sum_{ \substack{   n \in \smooth(X)} }
\frac{\alpha_{-1}(gn)}{n} \right)^2
+ O\left( d(g)^2 g^{1/2}  T^{-1/200} \right)\,.
\end{align}
By  arguments similar to those above we also find that
$$
\sum_{ \substack{ g \in \smooth(X)}} \frac{\phi(g) d(g)^2}{g^{3/2}}  \ll  T^{1/400}
\quad \hbox{and}\quad
\sum_{ \substack{g>T^{1/20} \\  g \in \smooth(X)}} \frac{\phi(g) d(g)^2}{g^2} \ll
 T^{-1/100} \,.
 $$
Using these and \eqref{sum squared}, we find that the right-hand side of
 \eqref {main term sum 1} equals
\begin{align*}
\left(   \sum_{ g \in \smooth(X) }  -   \sum_{ \substack{g > T^{1/20} \\  g \in \smooth(X) }} \right)
&\frac{\phi(g)}{g^2}  \left(    \sum_{  n \in \smooth(X)}  \frac{\alpha_{-1}(gn)}{n} \right)^2
+O\left(   T^{-1/200} \sum_{ \substack{g \leq T^{1/20} \\  g \in \smooth(X)}}
 \frac{\phi(g) d(g)^2}{g^{3/2}}  \right)  \\
&= \sum_{ g \in \smooth(X) }
 \frac{\phi(g)}{g^2}  \left(    \sum_{  n \in \smooth(X)}  \frac{\alpha_{-1}(gn)}{n} \right)^2
 + O\left(   T^{-1/400} \right) \,.
\end{align*}
Combining this with \eqref{int formula 2}, we now have
\begin{align}\label{m value}
\frac{1}{T} \int_T^{2T} \left|\zeta(\tfrac12+\i t)\right|^2 &\Bigg|
\sum_{\substack{n \leq T^{1/20} \\ n \in \smooth(X)}}
\frac{\alpha_{-1}(n)}{n^{1/2+\i t}} \Bigg|^2 \;\d t \\
&= \log T \sum_{ g \in \smooth(X) }
 \frac{\phi(g)}{g^2}  \left(    \sum_{  n \in \smooth(X)}  \frac{\alpha_{-1}(gn)}{n} \right)^2
 + O\left((\log X)^{10}\right) \,.  \notag
\end{align}
Since $\alpha_{-1}$ and $\phi$ are multiplicative
functions, we may  expand the entire sum into the Euler product
\begin{equation*}
\prod_{p\leq X} \left(\sum_{r} \sum_{j} \sum_{k}
\frac{\varphi(p^r) \alpha_{-1}(p^{j+r}) \alpha_{-1}(p^{k+r})}{p^{2r+j+k}}
\right)\,.
\end{equation*}
Recall that
$\alpha_{-1}(n) = \mu(n)$, the M\"obius function, if $n\in
\smooth(\sqrt{X})$; $\alpha_{-1}(p) = \mu(p)$ for all $p\leq X$; and
$\alpha_{-1}(n) \ll d_{3/2}(n) \ll d(n)$ for all $n\in \smooth(X)$.
Thus, the product equals
\begin{align*}
\prod_{p\leq \sqrt{X}} \left( 1 - \frac{1}{p} \right)
\prod_{\sqrt{X} < p \leq X} \left(1-\frac{1}{p} + O\left(\frac{1}{p^2}\right)\right)
&=\prod_{p\leq X} \left(1-\frac{1}{p}\right) \prod_{\sqrt{X}<p\leq
X} \left(1+O(\frac{1}{p^2})\right)\\
&=\frac{1}{e^{\gamma}\log X}\left(1 + O\left(\frac{1}{\log
X}\right)\right)\,.
\end{align*}
Since $\log X \ll \log\log T$, it now follows from \eqref{m value} that
\begin{align}\label{main term 3}
\frac{1}{T} \int_T^{2T} \left|\zeta(\tfrac12+\i t)\right|^2
\Bigg|\sum_{\substack{n \leq T^{1/20} \\ n \in \smooth(X)}}
\frac{\alpha_{-1}(n)}{n^{1/2+\i t}} \Bigg|^2 \;\d t
= \frac{\log T }{e^{\gamma}\log X}\left(1 + O\left(\frac{1}{\log
X}\right) \right)  \,.
\end{align}

Rewriting \eqref{eq:ranking splitting the sum 1} (with $\theta=1/20 $)
as $P_{X}^\ast(\tfrac12 + \i t)^{-1} =  \sum + \,O(T^{- \epsilon/200})$,
we see that
\begin{multline*}
\frac{1}{T}\int_T^{2T} |\zeta(\tfrac12 + \i t )
P_{X}^\ast(\tfrac12 + \i t)^{-1} |^{2 } \;\d t =
\frac{1}{T}\int_T^{2T} |\zeta(\tfrac12 + \i t ) |^{2 }
\left| \sum + O(T^{- \epsilon/200})\right|^2\;\d t  \\
=\frac{1}{T}\int_T^{2T} |\zeta(\tfrac12 + \i t ) |^{2 } \left|
\sum \right|^2\;\d t + O\left( \frac{1}{T^{1+\epsilon/200}}
\int_T^{2T} |\zeta(\tfrac12 + \i t ) |^{2 }
\left| \sum \right|\;\d t  \right)  \\
 + O \left( \frac{1}{T^{1+\epsilon/100}}
\int_T^{2T} | \zeta(\tfrac12 + \i t ) |^{2 }   \;\d t   \right) \,.
\end{multline*}
The final term is $O(T^{-\epsilon/200})$ since the second moment
of the zeta function is $O(T\log T)$. Also, by the Cauchy-Schwarz
inequality and \eqref{main term 3}, the second term is
\begin{align*}
\ll & \frac{1}{T^{1+\epsilon/200}}  \left(
\int_T^{2T} \left| \zeta(\tfrac12 + \i t )\;\sum \right|^{2 } \;\d t
\,\int_T^{2T} \left|\zeta(\tfrac12 + \i t ) \right|^{2 } \;\d t  \right)^{1/2}  \\
 \ll  & \frac{1}{T^{1+\epsilon/200}}  \left( T^2  \log^2 T /\log X  \right)^{1/2}
\ll T^{-\epsilon/400} \,.
\end{align*}
{}From these estimates and \eqref{main term 3}, we may now conclude that
\begin{align*}
\frac{1}{T}\int_T^{2T} |\zeta(\tfrac12 + \i t )
P_{X}^\ast(\tfrac12 + \i t)^{-1} |^{2 } \;\d t
=  \frac{\log T}{ e^{\gamma}\log X}\left(1+ O\left(\frac{1}{\log X}\right) \right)
\end{align*}
for $X=O(\log T)^{2-\epsilon}$.
This completes the proof of Theorem~\ref{thm:first case of splitting} in
the case $k=1$.


We now prove Theorem~\ref{thm:first case of splitting} for $k=2$.
By  Lemma~\ref{lem:restricting_P_to primes and squares} we may
again replace $ P_X(\tfrac12+\i t)$ by $ P_X^\ast(\tfrac12+\i t)$.
Furthermore, $G^2(3)/G(5)=1/12$, so it suffices to show that
\begin{equation}\label{k=2 splitting}
\frac{1}{T} \int_T^{2T} \left|
 \zeta(\tfrac12+\i t)^2  P_X^\ast(\tfrac12+\i t)^{-2}\right|^{2}
 \;\d t = \frac{1}{12} \left(1+ o\left(1\right) \right)
\left(  \frac{\log T}{e^{\gamma} \log X}\right)^4
\end{equation}
for $X \ll (\log T)^{2-\epsilon}$.
By \eqref{eq:tmp_truncate_smooth_sum_in_P}
 (see \eqref{eq:ranking splitting the sum 1} also and the remark following it), we have
\begin{equation*}\label{eq:ranking splitting the sum 2}
P_X^\ast(\tfrac12 +\i t)^{-2} = \sum_{\substack{n \leq T^{\theta} \\ n \in
\smooth(X)}} \frac{\alpha_{-2}(n)}{n^{1/2+\i t}} +
O\left(T^{- \epsilon\theta/10} \right) \,,
\end{equation*}
say, where $\alpha_{-2}(p) = -2$ for all $p\leq X$,
$\alpha_{-2}(p^2) = 1$ if $p \leq \sqrt X$,
$\alpha_{-2}(p^2) = 2$ if $ \sqrt X < p \leq X$,
and  $\alpha_{-2}(p^j) =0$ otherwise. In particular,
we note that $|\alpha_{-2}(n)| \leq d(n)$.

In carrying out the proof of splitting for this case, we  will gloss over some
of the less important steps as these are handled analogously to those for the
$k=1$ case. In particular, by an argument similar
to the one at the end of the proof of the case $k=1$, one can show that
\begin{align}\label{formula -2}
\frac{1}{T} \int_T^{2T} & \left|
 \zeta(\tfrac12+\i t)^2  P_X^\ast(\tfrac12+\i t)^{-2} \right|^{2}  \;\d t  \\
 &=
 \left(1+ O\left(\frac{1}{\log X}\right)\right)
 \frac{1}{T} \int_T^{2T} \left| \zeta(\tfrac12+\i t)^2
 \sum_{\substack{n \leq Y \\ n \in
\smooth(X)}} \frac{\alpha_{-2}(n)}{n^{1/2+\i t}}   \right|^{2} \;\d t \,, \notag
\end{align}
where $Y=T^\theta$ and $\theta>0$. Eventually we will take
$\theta$ very  small.

To estimate the right-hand side we use an analogue of \eqref{int formula} due to
Jose Gaggero~\cite{JG}. Let $A(s) = \sum_{n\leq Y} a_n n^{-s}$, where the $a_n$
are complex coefficients and $Y =T^\theta$ with $\theta < 1/150$.
Gaggero's  formula is
\begin{align*} 
\left( 1+O(\frac{1}{(\log T)^B}) \right) & \frac1T
 \int_{T}^{2T} \left| \zeta(\tfrac12 +\i t)\,A(\tfrac12 +\i t) \right|^2\,\d t     \\
  = &  \sum_{k=1}^{4}\left\{
\sum_{m, n \leq Y} \frac{c_k(m,n) a_m \overline{a_n }}{mn}
(m, n) \left(
\log^k\left( \frac{YT(m,n)}{2\pi mn}\right)
+\log^k \left(\frac{ T(m,n)}{2\pi Y}\right)  \right) \right\}  \notag  \\
 - &   \sum_{m, n \leq Y} \frac{ a_m \overline{a_n }}{mn}
\sum_{0< d<Y/4} \frac{(m,d)(n,d)}{d}   \left( \log(\frac{Y}{4d}) + O(1) \right)
\sum_{v<V_1} \frac1v  \sum_{\substack{u<U_1 \\ (n_d u, m_d v) =1}} \frac1u   \notag \\
 -  &  \sum_{m, n \leq Y} \frac{ a_m \overline{a_n }}{mn}
\sum_{0< d<mn/4Y} \frac{(m,d)(n,d)}{d} \left(  \log(\frac{mn}{4dY}) + O(1) \right)
\sum_{v<V^{\prime}_1} \frac1v
\sum_{\substack{u<U^{\prime}_1 \\ (n_d u, m_d v) =1} } \frac1u    \notag \,.
\end{align*}
Here
\begin{align}\label{U,V defns}
&U_1= CYT/  d n_d, \qquad\,\;V_1= CYT/  d m_d, \\
&U^{\prime}_1= CmnT/ Y d n_d, \quad V^{\prime}_1= CmnT/ Y d m_d \,,
\notag
\end{align}
 $C=2/\pi$, $B$ is an arbitrary positive number,
and for integers $n$ and $d$ we write $n_d=n/(n, d)$. Also,
$c_4(m, n) =(1/4\pi^2) \delta(m_n) \delta(n_m)$, where
\begin{equation*}
\delta(n)=\prod_{p^r \mid\mid n}\left(1+ r
\frac{(1-1/p)}{(1+1/p)}\right),
\end{equation*}
and $c_j(m, n) \ll |c_4(m, n)| (\log \log 3mn  )^{4-j}$ for $j=1,
2, 3$.

To estimate the right-hand side of \eqref{formula -2}, we take  $a_n=\a_{-2}(n)$
and $Y=T^{\epsilon_1}$ in this and obtain
\begin{align}\label{Jose's Formula 2}
 &\left(1+  O\left(\frac{1}{\log X}\right)\right)
 \frac{1}{T} \int_T^{2T}  \left|
 \zeta(\tfrac12+\i t)^2 P_X^\ast(\tfrac12+\i t)^{-2}  \right|^{2}  \;\d t   \notag\\
 &= \left(\frac{1}{2\pi^2} + O(\epsilon_1) \right)\log^4T \sum_{\substack{m, n \leq Y\\
 m, n \in \smooth(X)} }
\frac{  \a_{-2}(m) \a_{-2}(n) \delta(m/(m, n)) \delta(n/(m, n))}{mn} (m, n) \notag   \\
& \quad - \sum_{\substack{ m, n \leq Y \\ m, n \in \smooth(X)}}
\frac{ \a_{-2}(m) \a_{-2}(n)}{mn} \sum_{\substack{ 0< d<Y/4   }}
\frac{(m,d)(n,d)}{d} \left( \log(\frac{Y}{4d}) + O(1) \right)
\sum_{v<V_1} \frac1v  \sum_{\substack{u<U_1 \\ (n_d u, m_d v) =1}} \frac1u    \\
& \quad - \sum_{\substack{ m, n \leq Y \\ m, n \in \smooth(X)}}
\frac{ \a_{-2}(m) \a_{-2}(n)}{mn} \sum_{\substack{ 0< d< mn/4Y }}
\frac{(m,d)(n,d)}{d} \left( \log(\frac{mn}{4dY}) + O(1) \right)
\sum_{v<V^{\prime}_1} \frac1v  \sum_{\substack{u<U^{\prime}_1
 \\ (n_d u, m_d v) =1}} \frac1u    \notag \\
 &= \T_1 - \T_2 - \T_3,  \notag
\end{align}
say.

Let us denote the sum in $\T_1$ by $S_1$. Grouping together those
terms for which $(m,n)= g$ and then replacing
$m$ by $mg$ and $n$ by $ng$, we obtain
\begin{align}\label{S_1}
S_1 =
 \sum_{\substack{ g \leq Y \\  g \in \smooth(X)}} \frac{1}{g}
\sum_{\substack{ n \leq Y/g \\ n \in  \smooth(X)}} \frac{  \a_{-2}(gn)
\delta\left(n\right)}{n}
\left( \sum_{\substack{ m \leq Y/g \\ (m,n) =1 \\ m \in  \smooth(X) }}
\frac{  \a_{-2}(gm) \delta\left(m\right)}{m}   \right) \,.
\end{align}

Let $P = \prod_{p \leq X} p$\,.
Since $\a_{-2}$ is supported on cube-free integers, the $g$'s we are summing
over may be restricted to numbers of the form
$$
g= g_{1}\, g_{2}^2, \quad\hbox{where} \quad g_1 \mid P,  \quad g_2 \mid (P/g_1) \,.
$$
Note that this representation is unique and that $(g_1, g_2) =1$.
The summation over $g$ in \eqref{S_1} may therefore be replaced by the double sum
\begin{equation*}
  \sum_{\substack{ g_1 \leq Y \\  g_1 \mid P }}
    \sum_{\substack{ g_2 \leq \left( Y/ g_1 \right)^{\frac12}\\
    g_2 \mid (P/g_1) }} \,.
\end{equation*}

In the sum over $n$ we group terms together according to their greatest
common divisor with $g= g_{1}\, g_{2}^2$.
Observe that we may assume that $(n, g_2) = 1$, for otherwise a cube  divides
$  g_{1}\, g_{2}^2 \,n $  and   $\a_{-2}(gn)  $   vanishes.
If we then write $(n, g_1) = r$
and $n=r N$, we may replace the sum over $n$ in \eqref{S_1}
by
\begin{equation*}
\sum_{\substack{ r \mid g_1   }}
\sum_{\substack{ N \leq (Y/  r g_1 g_2^2 ) \\ N \in \smooth(X)\\ (N, (g_1/r) g_2)  =1}} \,.
\end{equation*}

Ignoring  the restriction $(m, n) =1$ for the moment, we may
similarly write the sum over $m$ in \eqref{S_1} as
\begin{equation*}\label{m sum}
\sum_{\substack{ s \mid g_1   }}
\sum_{\substack{ M \leq (Y/  s g_1 g_2^2 ) \\ M \in \smooth(X)\\ (M, (g_1/s) g_2)  =1}} \,.
\end{equation*}
Instead of  $(m, n) =1$ we now have $(sM, rN) =1$ or, equivalently,
$(M, N) = (r, s) = (N, s) = (M, r) =1$. We may impose the condition $(r, s)=1$
by replacing   $s \mid g_1$ in \eqref{m sum} by
$s \mid (g_1/r)$ since $g_1$ is square-free. Furthermore, since
$(N, g_1/r)=1$ and $s \mid (g_1/r)$, we automatically have  $(N, s)=1$.
 Thus, the coprimality conditions on $M$
are $(M, (g_1/s) g_2)  = (M, N) = (M, r)=1$. The first condition
implies the third because $r\mid (g_1/s)$. Thus, we need only
require that $(M, N(g_1/s) g_2) =1$.
The sum over $m$ may therefore be written
\begin{equation*}
\sum_{\substack{ s \mid g_1/r}}
\sum_{\substack{ M \leq (Y/  s g_1 g_2^2 ) \\ M \in \smooth(X)\\
(M, N(g_1/s) g_2)  =1}} \,.
\end{equation*}

We now have
\begin{align*}
S_1  =
 \sum_{\substack{ g_1 \leq Y \\  g_1 \mid P }} \frac{1}{g_1}
    \sum_{\substack{ g_2 \leq \left( Y/ g_1 \right)^{\frac12}\\
    g_2 \mid (P/g_1) }}  \frac{1}{g_{2}^2}
 \sum_{\substack{ r \mid g_1   }}  \frac{1}{r}
 \sum_{\substack{ s \mid (g_1/r)   }}    \frac{1}{s}
 & \sum_{\substack{ N \leq (Y/  r g_1 g_2^2 ) \\ N \in \smooth(X) \\
  (N, (g_1/r) g_2)  =1}}     \frac{  \a_{-2}(r^2 g_2^2 N (g_1/r))
 \delta\left(rN\right) }{N}  \\
 &\quad \times \sum_{\substack{ M \leq (Y/  s g_1 g_2^2 ) \\ M \in \smooth(X)\\
(M, N(g_1/s) g_2)  =1}}
\frac{\a_{-2}(s^2 g_2^2 M (g_1/s))
 \delta\left(s M\right)  }{M}    \,.\notag
\end{align*}
Note that if $N$ and $r$ have a common factor, then $\a_{-2}(r^2 g_2^2 N (g_1/r)) =0$,
and similarly for $M$ and  $s$. We may  therefore  replace the coprimality conditions in the sums over $N$ and $M$ by   $(N, g_1 g_2)  =1$
and $(M, N g_1 g_2)  =1$, respectively. The new conditions then imply that
$\a_{-2}(r^2 g_2^2 N (g_1/r)) = \a_{-2}(r^2)
\a_{-2}(g_2^2) \a_{-2}(N) \a_{-2}(g_1/r) ,  \delta(rN) = \delta(r)\delta(N)$,
and similarly for  $\a_{-2}(s^2 g_2^2 M (g_1/s))$ and $\delta(sM)$.  Hence
\begin{align*}
S_1 &=
 \sum_{\substack{ g_1 \leq Y \\  g_1 \mid P }} \frac{\a_{-2}(g_1)^2 }{g_1}
\sum_{\substack{ g_2 \leq \left( Y/ g_1 \right)^{\frac12}\\
g_2 \mid (P/g_1) }}  \frac{\a_{-2}(g_2^2)^2}{g_{2}^2}
 \sum_{\substack{ r \mid g_1   }}  \frac{ \a_{-2}(r^2) \delta(r)}{ \a_{-2}(r)r}
 \sum_{\substack{ s \mid (g_1/r)   }}    \frac{\a_{-2}(s^2) \delta(s)}{\a_{-2}(s)s}  \\
& \qquad \sum_{\substack{ N \leq (Y/  r g_1 g_2^2 ) \\ N \in \smooth(X)
\\ (N, g_1 g_2)  =1}}
 \frac{\a_{-2}(N) \delta(N)}{N}
\sum_{\substack{ M \leq (Y/  s g_1 g_2^2 ) \\ M \in \smooth(X)\\
(M, N g_1 g_2)  =1}}
 \frac{\a_{-2}(M) \delta(M)}{M}   \,.  \notag
\end{align*}

We next extend each of the sums here to all of $\smooth(X)$.
The error terms this introduces are handled
as  they were in the case $k=1$,  and  they contribute
at most  ``little $o$'' of the main term. Observing also that
$M$ and $N$ may be restricted to cube-free integers, we obtain
\begin{align}\label{S_1*}
S_1  = \left( 1+ o(1)  \right) &
 \sum_{\substack{  g_1 \mid P }} \frac{\a_{-2}(g_1)^2 }{g_1}
\sum_{\substack{g_2 \mid (P/g_1) }}  \frac{\a_{-2}(g_2^2)^2}{g_{2}^2}
 \sum_{\substack{ r \mid g_1   }}  \frac{ \a_{-2}(r^2) \delta(r)}{ \a_{-2}(r)r}
 \sum_{\substack{ s \mid (g_1/r)   }}    \frac{\a_{-2}(s^2) \delta(s)}{\a_{-2}(s)s}  \\
& \qquad \sum_{\substack{ N \mid (P/g_1g_2)^2}}
 \frac{\a_{-2}(N) \delta(N)}{N}
\sum_{\substack{   M \mid (P/N g_1 g_2)^2}}
 \frac{\a_{-2}(M) \delta(M)}{M}   \,.  \notag
\end{align}

We now define the following multiplicative functions:
\begin{align*}
A(n) &= \sum_{\substack{ d\mid n}}
 \frac{\a_{-2}(d) \delta(d)}{d}=
 \prod_{p^a\mid\mid n}\left(1 + \frac{ \a_{-2}(p) \delta(p) }{p }
 \cdots +\frac{ \a_{-2}(p^a) \delta(p^a) }{p^a} \right),\\
B(n) &= \sum_{\substack{ d\mid n}}
 \frac{\a_{-2}(d) \delta(d)}{d A(d^2)}=
\prod_{p^a\mid\mid n}\left(1 + \frac{ \a_{-2}(p) \delta(p) }{p
\,A(p^2)} + \cdots + \frac{ \a_{-2}(p^a) \delta(p^a) }{p^a
\,A(p^{2a})}  \right),\\
 C(n) &=    \sum_{\substack{ d \mid (n, P)
}}    \frac{\a_{-2}(d^2) \delta(d)}{\a_{-2}(d)d} =\prod_{p\mid (n,
P)}\left(1 + \frac{ \a_{-2}(p^2) \delta(p) }{\a_{-2}(p) p}
\right),\\
 D(n) &=   \sum_{\substack{ d \mid (n, P) }}
\frac{\a_{-2}(d^2) \delta(d)}{\a_{-2}(d) C(d)d} =\prod_{p\mid (n,
P)}\left(1 + \frac{ \a_{-2}(p^2) \delta(p) }{\a_{-2}(p)C(p) p}
\right),\\
E(n) &=  \sum_{d\mid n} \frac{ \a_{-2}(d^2 )^2 } {  A(d^2) B(d^2
)\,d^2 } = \prod_{p\mid l}\left(1 + \frac{ \a_{-2}(p^2 )^2 }
{A(p^2) B(p^2 )\,p^2 } \right),
\end{align*}
and
\begin{equation*}
F(n) = \sum_{d\mid n} \frac{ \a_{-2}(d)^2 C(d)D(d)} { A(d^2) B(d^2
) E(d)\,d} = \prod_{p\mid l}\left(1 + \frac{ \a_{-2}(p )^2
C(p)D(p)} { A(p^2) B(p^2 ) E(p)\,d} \right).
\end{equation*}

Using these definitions and working from the inside out in \eqref{S_1*},
we find first that the sum over $M$ is $A((P/Ng_1g_2)^2) =
A(P^2)/A(N^2)A(g_1^2)A(g_2^2)$. The contribution of the sums
over $M$ and $N$ together is then $\left(A(P^2)/A(g_1^2)A(g_2^2)\right)
\left( B(P^2)/ B(g_1^2)B(g_2^2)\right)$. Thus, so far we have
\begin{align*}
S_1  = \left( 1+ o(1)  \right)    A(P^2)B(P^2)
 \sum_{\substack{  g_1 \mid P }} &\frac{\a_{-2}(g_1)^2 }{g_1A(g_1^2)B(g_1^2)}
\sum_{\substack{g_2 \mid (P/g_1) }}
\frac{\a_{-2}(g_2^2)^2}{g_{2}^2 A(g_2^2)B(g_2^2)} \\
& \sum_{\substack{ r \mid g_1   }}  \frac{ \a_{-2}(r^2) \delta(r)}{ \a_{-2}(r)r}
  \sum_{\substack{ s \mid (g_1/r)   }}    \frac{\a_{-2}(s^2) \delta(s)}{\a_{-2}(s)s}
   \,.  \notag
\end{align*}
The sums over $r$ and $s$ contribute $C(g_1) D(g_1)$,
and  the sum over $g_2$ is then $E(P)/E(g_1)$. Thus, we see that
\begin{align*}
S_1 &= \left( 1+ o(1)  \right)  A(P^2)B(P^2)E(P)
 \sum_{\substack{  g_1 \mid P }} \frac{\a_{-2}(g_1)^2
 C(g_1) D(g_1)}{g_1A(g_1^2)B(g_1^2) E(g_1)} \\
 &= \left( 1+ o(1)  \right)   A(P^2)B(P^2)E(P)F(P) \,.  \notag
\end{align*}
Using the expression for $F(P)$ as a product, we  see that this is the same  as
\begin{equation}\label{S_1  formula 1}
S_1 = \left( 1+ o(1)  \right)   \prod_{p\mid P}
\left( A(p^2)B(p^2)E(p)  + \frac{\a_{-2}(p)^2 C(p)D(p)}{p} \right)  \,.
\end{equation}
By the  definitions of $C$ and $D$ we see that
\begin{align}\label{CD}
C(p)D(p)
& =  \left(1 + \frac{ \a_{-2}(p^2) \delta(p) }{\a_{-2}(p) p}  \right)
+ \frac{ \a_{-2}(p^2) \delta(p) }{\a_{-2}(p) p}  \\
 &= 1 - \frac{ \a_{-2} (p^2) \delta(p) }{p} \,, \notag
\end{align}
since $ \a_{-2}(p)=-2$ for $p$ dividing $P$.
Similarly,
\begin{align*}
A(p^2)\,B(p^2)\,E(p)
&= A(p^2)\,B(p^2) + \frac{ \a_{-2}(p^2 )^2 }{p^2 } \,.
\end{align*}
It is clear that
$A(p^2)=A(p^3)= \cdots$. Therefore
\begin{align*}
B(p^2) &= 1 + \frac{ \a_{-2}(p) \delta(p) }{p \,A(p^2)}
+ \frac{ \a_{-2}(p^2) \delta(p^2) }{p^2 \,A(p^{2})}  \\
&= 1 + \frac{1}{A(p^{2})}\left( A(p^{2}) -1  \right) \\
&= 2 - \frac{1}{A(p^{2})}
\end{align*}
and
\begin{align*}
A(p^2)\,B(p^2)\,E(p)
&= 2A(p^2)-1 + \frac{ \a_{-2}(p^2 )^2 }{p^2 }  \,.
\end{align*}
We use this, \eqref{CD}, and $\a_{-2}(p)= -2$, and obtain
\begin{align*}
A(p^2)&B(p^2) E(p) + \frac{\a_{-2}(p)^2 C(p)D(p)}{p} \\
&=  2A(p^2) -1 + \frac{ \a_{-2}(p^2 )^2 }{p^2 }
+ \frac{ 4 }{p}  - \frac{  4 \a_{-2}(p^2 )^2\delta(p) }{p^2 }  \\
&= 2\left( 1 - \frac{ 2 \delta(p) }{p  }
+ \frac{ \a_{-2}(p^2) \delta(p^2) }{p^2}   \right)
 -1 + \frac{ \a_{-2}(p^2 )^2 }{p^2 }
+ \frac{ 4 }{p}  - \frac{ 4 \a_{-2}(p^2 )^2\delta(p) }{p^2 } \\
&= 1 + \frac{ 4-4 \delta(p) }{p  }
+ \frac{ \a_{-2}(p^2) \left(\a_{-2}(p^2) -4 \a_{-2}(p^2)\delta(p)
+2\delta(p^2) \right) }{p^2} \,.
\end{align*}
Recall that $\delta(p^r) = 1+ r \frac{(1-1/p)}{(1+1/p)}$, so that
$\delta(p)= 2/(1+1/p)$ and $\delta(p^2)= 2\delta(p)-1$. Also recall that
$\a_{-2}(p^2)=1$ if $p \leq \sqrt{X}$.
Thus, for $p \leq \sqrt{X}$ the last line is
\begin{align*}
&=1 + \frac{ 4-4 \delta(p) }{p  }
+ \frac{   1 -4  \delta(p) + 2 \left(2\delta(p)-1 \right) }{p^2}
=1 + \frac{ 4-4 \delta(p) }{p  }
- \frac{ 1 }{p^2} \\
&=1 + \frac{ 4}{p  } -  \frac{ 8}{p+1} - \frac{ 1}{p^2}
=\frac{(1-1/p)^3}{1+1/p}
=\frac{(1-1/p)^4}{1-1/p^2}\,.
\end{align*}
On the other hand, if  $\sqrt{X} < p \leq X$, then  $\a_{-2}(p^2)=2$,
and the  last line  is
$$
=\frac{(1-1/p)^4}{1-1/p^2} + O\left(  1/p^2   \right)\,.
$$
Combining these results in  \eqref{S_1  formula 1},
we find that
\begin{align*}
S_1 &= \left( 1+ o(1)  \right)
\prod_{p\leq \sqrt{X}}\left( \frac{(1-1/p)^4}{1-1/p^2} \right)
\; \prod_{\sqrt{X}< p\leq X}\left( \frac{(1-1/p)^4}{1-1/p^2}
+ O\left(  1/p^2   \right)\right) \\
&= \left( 1+ o(1)  \right)
\prod_{p\leq  X}\left( \frac{(1-1/p)^4}{1-1/p^2} \right)
\; \prod_{\sqrt{X}< p\leq X}\left( 1+ O\left(  1/p^2   \right)\right) \\
&= \left( 1+ o(1)  \right)
\prod_{p\leq  X}\left( 1-1/p\right)^4
\prod_{p}\left(1-1/p^2\right)^{-1} \\
&= \left( 1+ o(1)  \right)  \frac{\pi^2}{6}\left(e^\gamma \log X \right)^{-4} \,.
\end{align*}
Since
$$
\T_1 =  S_1\, \left(\frac{1}{2\pi^2} + O(\epsilon_1) \right)\log^4T
$$
and $\epsilon_1>0$ may be taken as small as we like, we now see that
\begin{equation}\label{T_1 estimate}
\T_1 = \left(\frac{1}{12 } + o(1) \right)  \left(\frac{\log T}{e^\gamma \log X}\right)^4 \,.
\end{equation}


To treat  the second term on the right-hand side of
\eqref{Jose's Formula 2}, $\T_2$,
 we require  two lemmas.

\begin{lem}\label{Sum 1}
Suppose that $a$ and $b$ are positive integers with $(a, b)=1$. Then
for $b \leq x$, we have
\begin{equation*}
\sum_{\substack{n\leq x \\ (an, b)=1}} \frac{1}{n} =
\frac{\phi(b)}{b}\log x + O(\log \log 2b)\,.
\end{equation*}
\end{lem}
\begin{proof} Since $(a, b)=1$, the condition $(an, b)=1$ is  equivalent to
$(n, b)=1$. Thus, the sum is
\begin{align*}
\sum_{\substack{n\leq x \\ (n, b)=1}} \frac{1}{n}
= \sum_{n\leq x}\frac{1}{n}  \sum_{\substack{d|n\\d|b}} \mu(d)
=\sum_{d|b} \frac{\mu(d)}{d}\sum_{m\leq x/d} \frac1m
=\sum_{d|b} \frac{\mu(d)}{d} \left(\log x/d + O(1)\right) \,.
\end{align*}
Now
$\sum_{d|b}  \mu(d)/d= \phi(b)/b $
and
 $$
 \sum_{d|b}  \frac{\mu(d)\log d}{d}
 = \frac{\phi(b)}{b}\sum_{p|b}\frac{\log p}{p-1}
 \ll  \frac{\phi(b)}{b}\log \log 2b
 \ll  \log \log 2b\,.
 $$
 Furthermore,
 $$
 \sum_{d|b} \frac{|\mu(d)|}{d}
 = \prod_{p\mid b}\left(1+\frac{1}{p}\right)
 \leq \frac{b}{\phi(b)} \ll \log\log 2b \,.
 $$
 Thus we find that
\begin{align*}
\sum_{\substack{n\leq x \\ (an, b)=1}} \frac{1}{n}
= \frac{\phi(b)}{b}\log x  + O(\log \log 2b)  \,.
\end{align*}
\end{proof}

\begin{lem}\label{Sum 2}
Let $\kappa(n)= \prod_{p\mid n}\left(  1+\frac1p \right)^{-1}$ and let $U, V$ be
either $U_1, V_1$ or $U^{\prime}_1, V^{\prime}_1$ as defined in
\eqref{U,V defns}. If $m, n, d \ll Y\leq T^{1/150} $, and $(n_d , m_d ) =1$,
then
\begin{equation*}
\sum_{v<V} \frac1v  \sum_{\substack{u<U \\ (n_d u, m_d v) =1} } \frac1u
= \frac{6}{\pi^2}\, \kappa(m_d)\,\kappa(n_d) \log U \,\log V + O(\log T \log\log T)\,.
\end{equation*}
\end{lem}
\begin{proof}  The conditions $(n_d u , m_d v ) =1$ and $(n_d, m_d) =1$ are equivalent to
$(v, n_d)=1$ and $(u, m_d v)=1$. Hence, by Lemma~\ref{Sum 1},
the double sum equals
\begin{align}\label{MT sum}
\sum_{\substack{v<V\\(v, n_d)=1} }\frac1v
\sum_{\substack{u<U \\ ( u, m_d v) =1} } \frac1u
&=  \log U  \sum_{\substack{v<V\\(v, n_d)=1} }\frac1v
\left( \frac{\phi(m_d v)}{m_d v}
+ O\left(\log \log (m_d V)\right) \right) \notag \\
&= \log U \sum_{\substack{v<V\\(v, n_d)=1} }\frac1v
\left( \frac{\phi(m_d v)}{m_d v} \right)
+   O\left( \log T \log\log T\right) \,.
\end{align}
Denoting the sum on the right by $\sum$, we have
\begin{align*}
\sum= \sum_{\substack{v<V\\(v, n_d)=1} }\frac1v
  \sum_{r| m_d v}\frac{\mu(r)}{r}
=\sum_{r<m_d V}\frac{\mu(r)}{r}    \sum_{\substack{v<V\\(v, n_d)=1\\ r| m_d v} }
\frac{1}{v} \,.
\end{align*}
Now set $(m_d, r)=g$ and write $r=gR$. Then $(m_d /g, R)=1$ and we find that
\begin{align*}
\sum=
 \sum_{g|m_d} \frac{1}{g}\sum_{\substack{R<m_d V/g\\ (m_d/g,R) = 1}}
\frac{\mu(gR)}{R}   \sum_{\substack{v<V\\(v, n_d)=1\\ R | v} }
\frac{1}{v}  \,.
\end{align*}
If we set $v=Rw$, then $w<V/R$, and
$(Rw, n_d)=1$ is the same as the two conditions $(R, n_d)=1$ and  $(w, n_d)=1$.
Thus, using Lemma~\ref{Sum 1} and the observation that the inner sum vanishes
unless $R<V$, we obtain
\begin{align*}
 \sum &= \sum_{g|m_d} \frac{1}{g}\sum_{\substack{R<m_d V/g\\ (m_d/g, R) = 1 \\ (n_d, R) = 1}}
\frac{\mu(gR)}{R^2}  \sum_{\substack{w<V/R\\(w, n_d)=1} }
\frac{1}{w}   \notag  \\
&= \sum_{g|m_d} \frac{1}{g}\sum_{\substack{R<  V \\ (m_d/g, R) = 1 \\ (n_d, R) = 1}}
\frac{\mu(gR)}{R^2} \left( \frac{\phi(n_d)}{n_d} \log \frac{V}{R } + O(\log\log 2 n_d)  \right)\,.
\end{align*}
We may assume $(R, g)=1$, for otherwise $\mu(gR)=0$. The coprimality conditions on the sum
may then be written $(m_d n_d, R)=1$, and we find that
\begin{align*}
 \sum &= \sum_{g|m_d}\frac{\mu(g )}{g}
  \sum_{\substack{R<  V \\ (R, m_d n_d) = 1}}
\frac{\mu(R)}{R^2} \left( \frac{\phi(n_d)}{n_d} \log \frac{V}{R}  + O(\log\log 2 n_d)  \right)
\notag \\
&= \frac{\phi(n_d)}{n_d} \log V   \sum_{g|m_d}\frac{\mu(g )}{g}
  \sum_{\substack{R<  V \\ (R, m_d n_d) = 1}}
\frac{\mu(R)}{R^2}
+ O\left( \log\log 2 n_d   \sum_{g|m_d}\frac{|\mu(g )|}{g}
  \sum_{ R<  V} \frac{\log R}{R^2}\right).
\end{align*}
Since $\sum_{g|m_d} |\mu(g )|/g =\prod_{p|m_d}(1+1/p) \ll \log\log 2 m_d$,
the error term is $\ll \left(\log\log  2 m_d \,  \log\log  2 n_d \right)$. The main term is
\begin{align*}
&= \frac{\phi(n_d)}{n_d} \log V   \sum_{g|m_d}\frac{\mu(g )}{g}
  \left( \sum_{\substack{R=1 \\ (R, m_d n_d) = 1}}^{\infty}
\frac{\mu(R)}{R^2} +  O( V^{-1}) \right)   \notag \\
&=   \zeta(2)^{-1} \prod_{p|m_d n_d}\left(1-\frac{1}{p^2}\right)^{-1}
\frac{\phi(n_d)}{n_d} \log V   \sum_{g|m_d}\frac{\mu(g )}{g}
+ O\left( \frac{\log V}{V}  \sum_{g|m_d}\frac{|\mu(g )|}{g}  \right) \\
&= \frac{6}{\pi^2}  \prod_{p|m_d n_d}\left(1-\frac{1}{p^2}\right) ^{-1}
\frac{\phi(m_d)}{m_d}\frac{\phi(n_d)}{n_d} \log V
+ O\left( \frac{\log V \log\log 2m_d }{V}    \right)     \notag  \,.
\end{align*}
By hypothesis, $(m_d, n_d)=1$. Furthermore,
$ \prod_{p| l}\left(1-1/p^2\right) ^{-1} (\phi(l)/l)
=  \prod_{p| l}\left(1+ 1/p\right) ^{-1} =\kappa(l) $.
Thus, combining our estimates, we obtain
$$
\sum= \frac{6}{\pi^2} \kappa(m_d) \kappa(n_d) \log V + O (\log\log 2 m_d  \log\log 2 n_d)\,.
$$
Since $m, n \ll T^{1/150}$, we obtain from this and
 \eqref{MT sum} that
$$
\sum_{v<V} \frac1v  \sum_{\substack{u<U \\ (n_d u, m_d v) =1} } \frac1u
=\frac{6}{\pi^2} \kappa(m_d) \kappa(n_d) \log U \log V + O (\log T \log\log T)\,.
$$
\end{proof}

Returning  to  $\T_2$ in \eqref{Jose's Formula 2} and using
Lemma~\ref{Sum 2}, we have
\begin{align*}
\T_2 =  \frac{6}{\pi^2}\, \sum_{\substack{ m, n \leq Y \\ m, n \in \smooth(X)}}
& \frac{  \a_{-2}(m) \a_{-2}(n)}{mn}
 \sum_{\substack{ 0< d<Y/4 \\  d\in \smooth(X) }}  \frac{(m,d)(n,d)}{d}   \left( \log(\frac{Y}{4d}) + O(1) \right) \\
&\qquad \times \bigg( \kappa(m_d)\,\kappa(n_d) \log U_1 \,\log V_1 + O(\log T \log\log T)
 \bigg)\,,
\end{align*}
where $U_1= CYT/  d n_d,  V_1= CYT/  d m_d,$ and $Y =T^\epsilon_1$.
Interchanging the order of  summation, we find that
\begin{align*}
\T_2 =  \frac{6}{\pi^2}\, \sum_{\substack{ 0< d<Y/4 \\  d \in \smooth(X) }}
\frac{1}{d} \left( \log(\frac{Y}{4d}) + O(1) \right)
&  \sum_{\substack{ m, n \leq Y \\ m, n \in \smooth(X)}}    \frac{  \a_{-2}(m) \a_{-2}(n)(m,d)(n,d) }{mn}  \\
 \times  & \bigg( \kappa(m_d)\,\kappa(n_d) \log U_1 \,\log V_1
 + O(\log T \log\log T) \bigg)\,.
\end{align*}
Since $\kappa(n) \ll \log\log 3n$, the  expression in the last parentheses is
$$
=\kappa(m_d)\,\kappa(n_d) \left( \log U_1 \,\log V_1   +
O(\log T \log\log^3 T)\right) =
(1+O(\epsilon_1))\kappa(m_d)\,\kappa(n_d)\, \log^2 T \,.
$$
Thus,
\begin{align}\label{T_2}
\T_2
&=  \frac{6}{\pi^2}\, ( 1 + O(\epsilon_1) ) \log^2 T
\sum_{\substack{ 0< d<Y/4 \\  d \in \smooth(X) }}
\frac{ (\log(Y/4d) + O(1) ) }{d}
\left( \sum_{\substack{  n \leq Y \\   n \in \smooth(X)}}  \frac{\a_{-2}(n)\,(n,d) \kappa(n/(n,d)) }{n}
  \right)^2  .
\end{align}

Denote the inner sum by $S(d)$. As on previous occasions, extending the sum
to all of $\smooth(X)$, we introduce an error term that is $o(1)$ times the main term.
Thus,  grouping together  terms in $S(d)$ for which $(n, d) =e$, say, we obtain
\begin{align*}
S(d)
&= \left( 1 + o(1)  \right)\sum_{e\mid d} e
\sum_{ \substack{     n \in \smooth(X) \\ (n, d)=e}}
\frac{\a_{-2}(n)\,  \kappa(n/e) }{n}
 = \left( 1 + o(1)  \right) \sum_{e\mid d}
\sum_{ \substack{     N \in \smooth(X) \\ (N, d/e)=1}}
\frac{\a_{-2}(eN)\,  \kappa(N) }{N} \,.
\end{align*}
Since $\a_{-2}$ is supported only on cube-free numbers in
$\smooth(X)$, we may assume that
$e \mid P^2$. Therefore, $e\mid (d, P^2)=D$, say. Now $D$ may be written uniquely as
$D=D_1D_2^2,$ where $D_1\mid P$ and $D_2\mid (P/D_1)$, so that, in particular,
$(D_1, D_2)=1$. Furthermore, we may write any divisor $e$ of $D$ as
 $e = e_1 e_2 e_3^2$, where $e_1\mid D_1, e_2\mid D_2,$ and $e_3\mid (D_2/e_2)$.
 Note that this means the $e_i$ are  pairwise coprime. The condition $(N, d/e)=1$ is now
 $(N, (D_1D_2^2/e_1 e_2 e_3^2  ))=1$.  Also, $\a_{-2}(eN)= \a_{-2}(e_1 e_2 e_3^2 N)$,
so we may  assume that $(N, e_3)=1$ and, therefore, that
$(N, (D_1D_2^2/e_1 e_2 ))=1$.
Observe, moreover, that   $e_2\mid D_2$ implies $e_2\mid (D_2^2/e_2)$. Thus,
$(N, (D_1D_2^2/e_1 e_2 ))=1$ is the same as $(N, (D_1D_2 /e_1))=1$.
It follows that $N$ and $e_1$ can have a common factor,
 but not $N$ and $e_2$ or $e_3$. We may therefore write
 $\a_{-2}(e_1 e_2 e_3^2 N) =  \a_{-2}(e_1 N) \a_{-2}(e_2) \a_{-2}(e_3^2)$
 and
\begin{align*}
S(d) =S(D)
&= \left( 1 + o(1)  \right) \sum_{e_1\mid D_1} \sum_{e_2\mid D_2}  \a_{-2}(e_2)
\sum_{e_3\mid (D_2/e_2)} \a_{-2}(e_3^2)
\sum_{ \substack{     N \in \smooth(X) \\ (N, (D_1D_2 /e_1))=1}}
\frac{\a_{-2}(e_1 N)\,  \kappa(N) }{N} \\
&= \left( 1 + o(1)  \right) \sum_{e_1\mid D_1}
\sum_{ \substack{ N \in \smooth(X) \\ (N, (D_1D_2 /e_1))=1}}
\frac{\a_{-2}(e_1 N)\,  \kappa(N) }{N}
\sum_{e_3\mid D_2} \a_{-2}(e_3^2)
\sum_{e_2\mid (D_2/e_3)}  \a_{-2}(e_2) \,.
\end{align*}
The innermost sum  is
\begin{align*}
\sum_{e_2\mid (D_2/e_3)}  \a_{-2}(e_2)
&=\prod_{p\mid (D_2/e_3)}\left( 1+ \a_{-2}(p)\right)
 =\prod_{p\mid (D_2/e_3)}\left( 1-2\right) \\
&=\mu(D_2/e_3)=\mu(D_2) \mu(e_3).
\end{align*}
We also have
\begin{align*}
\sum_{e_3\mid D_2} \mu(e_3) \a_{-2}(e_3^2)
= \prod_{p \mid D_2}\left( 1- \a_{-2}(p^2)\right) \,.
\end{align*}
At this point it is convenient to define numbers
$$
P_1= \prod_{p\leq \sqrt{X}} p  \quad \hbox{and} \quad
P_2= \prod_{\sqrt{X}< p\leq X} p \,.
$$
Notice that $P=P_1 P_2$.
Since $\a_{-2}(p^2)=1$ if $p\mid P_1$ and  $\a_{-2}(p^2)=2$
if $p\mid P_2$, the sum over $e_3$ equals $0$ unless $D_2 \mid P_2$,
in which case it
equals $\mu(D_2)$. Thus, if $D_2$ and $P_1$ have a common factor,
$S(D_1D_2^2)=0$,
whereas if $D_2 \mid P_2$, then
 \begin{align*}
S(d)= S(D_1D_2^2)
&= \left( 1 + o(1)  \right) \sum_{e_1\mid D_1}
\sum_{ \substack{N \in \smooth(X) \\ (N, (D_1D_2 /e_1))=1}}
\frac{\a_{-2}(e_1 N)\,  \kappa(N) }{N} \,.
\end{align*}
From this point on we shall therefore assume that  $D_2 \mid P_2$.

Now set $(N, e_1)=r$ and write $N=rM$.
Then we have
\begin{align*}
S(D_1D_2^2)
&= \left( 1 + o(1)  \right) \sum_{e_1\mid D_1} \sum_{r \mid e_1}
\sum_{ \substack{ N \in \smooth(X) \\ (N,  e_1)= r  \\ (N, (D_1D_2 /e_1))=1}}
\frac{\a_{-2}(e_1 N)\,  \kappa(N) }{N} \\
&= \left( 1 + o(1)  \right) \sum_{e_1\mid D_1} \sum_{r \mid e_1}\frac{1}{r}
\sum_{ \substack{ M \in \smooth(X) \\ (M,  e_1/r)= 1  \\ (r M, (D_1D_2 /e_1))=1}}
\frac{\a_{-2}(r^2  M(e_1/r) )\,  \kappa(r M) }{M} \,.
\end{align*}
We may assume that $( M, r)=1$ and $(r, e_1/r) =1$, since otherwise
$\a_{-2}(r^2  M(e_1/r) )=0$. Actually, $(r, e_1/r) =1$
 is automatically satisfied because $r\mid e_1$ and $e_1$ is square-free.
It follows that  $\kappa(r M)=\kappa(r )\kappa(M)$ and, since we also have
$(M,  e_1/r)= 1$,   that
$\a_{-2}(r^2  M(e_1/r) )= \a_{-2}(r^2)  \a_{-2}(M) \a_{-2}(e_1/r)$.
The  coprimality conditions in the sum are now seen to be equivalent to
the conditions
$(M, r)=(r, e_1/r)=(M,  e_1/r)=( M, (D_1D_2 /e_1))
=(r, (D_1D_2 /e_1))=1$. As we have already pointed out,
the second of these is automatic.
Similarly, so is the last.
The remaining conditions are
equivalent to $( M, D_1D_2 )=1$, so we find that
\begin{align*}
S(D_1D_2^2) &= \left( 1 + o(1)  \right) \sum_{e_1\mid D_1} \a_{-2}(e_1)
\sum_{r \mid e_1}  \frac{\a_{-2}(r^2) \kappa(r) }{r \a_{-2}(r) }
\sum_{ \substack{ M \in \smooth(X) \\   ( M,  D_1D_2 )=1}}
\frac{\a_{-2}(M)\,\kappa(M)}{M} \,.
\end{align*}
The sum over $M$ equals
\begin{align}\label{G formula}
\prod_{p\mid (P/D_1D_2)} \left( 1 + \frac{\a_{-2}(p)\,\kappa(p)}{p}
+ \frac{\a_{-2}(p^2)\,\kappa(p^2)}{p^2} \right)
= G(P/D_1D_2)\,,
\end{align}
say. Hence,
\begin{align*}
S(D_1D_2^2) &= \left( 1 + o(1)  \right) G(P/D_1D_2)
\sum_{e_1\mid D_1} \a_{-2}(e_1)
\sum_{r \mid e_1}  \frac{\a_{-2}(r^2) \kappa(r) }{r \a_{-2}(r) }\,.
\end{align*}
The  double sum  equals
\begin{align} \label{H formula}
 \sum_{r \mid D_1} \frac{\a_{-2}(r^2) \kappa(r) }{r \a_{-2}(r) }
\sum_{f_1 \mid (D_1/r)} \a_{-2}(f_1)
&= \mu(D_1)  \sum_{r \mid D_1}
\frac{\mu(r) \a_{-2}(r^2) \kappa(r) }{r \a_{-2}(r) }  \notag \\
&=\mu(D_1)  \prod_{p \mid D_1} \left( 1   +
\frac{\a_{-2}(p^2) \kappa(p)}{2p} \right) \\
&= \mu(D_1) H(D_1)\,, \notag
\end{align}
say. Thus,
\begin{align}\label{S(d)}
S(d) = S(D_1D_2^2)  = \left( 1 + o(1)  \right) G(P)\frac{\mu(D_1) H(D_1)  }{G(D_1)G(D_2)} \,,
\end{align}
provided $D_2\mid P_2$; otherwise $S(d)=0$.

We  use this in  \eqref{T_2}. Recall that for each $d <Y/4$ we had set
$(d, P^2) =D_1D_2^2$ with $D_1\mid P$ and $D_2\mid (P/D_1)$. Recall
also that $P=P_1P_2$ and $Y=T^{\epsilon_1}$. We therefore  have that
\begin{align*}
\T_2
&=  \frac{6}{\pi^2}\, ( 1 + O(\epsilon_1) )  \log^2 T
\sum_{\substack{ 0< d<Y/4 \\  d \in \smooth(X) }}
\frac{ (\log(Y/4d) + O(1) ) }{d}  S(d)^2 \\
&= \frac{6}{\pi^2}\, ( 1 + O(\epsilon_1) )  \log^2 T
\sum_{D_2\mid P_2}   \frac{1}{D_2^2}\sum_{D_1\mid (P/D_2)} \frac{S(D_1D_2^2)^2 }{D_1}
\sum_{\substack{ 0< \delta<Y/4D_1D_2^2 \\  (\delta,  (P_1P_2)^2/D_1D_2^2 )=1 }}
\frac{ \left(\log(Y/(4D_1D_2^2 \delta)\right) + O(1) ) }{\delta} \,.
\end{align*}
The coprimality condition in the last sum is equivalent to
$(\delta,   P_1P_2 / D_2  )=1$.  Thus,
using \eqref{S(d)}, we find that
\begin{align*}
\T_2
 =  \frac{6}{\pi^2}\, ( 1 + O(\epsilon_1) )  G(P)^2 \log^2 T
\sum_{D_2\mid P_2} &  \frac{1}{D_2^2 G(D_2)^2}
\sum_{D_1\mid (P/D_2)} \frac{H(D_1)^2 }{D_1 G(D_1)^2} \\
& \times\sum_{\substack{ 0< \delta<Y/4D_1D_2^2 \\  (\delta,  P_1P_2 / D_2 )=1 }}
\frac{ \left( \log(Y/(4D_1D_2^2 \delta)) + O(1) \right) }{\delta}  \,.
\end{align*}
By Lemma~\ref{Sum 1}  the sum over $\delta$ is
 \begin{equation*}\label{delta sum}
\ll \log Y \sum_{\substack{ 0< \delta<Y/4D_1D_2^2 \\  (\delta,  P_1P_2 / D_2 )=1 }}
\frac{ 1 }{\delta}
 \ll      \frac{\phi(P) }{P} \frac{D_2}{\phi(D_2)}  \log^2 Y    \,.
\end{equation*}
Thus
\begin{align*}
\T_2
 \ll    G(P)^2 \frac{\phi(P) }{P} \log^2 T \log^2 Y
\sum_{D_2\mid P_2}  & \frac{1}{D_2  \phi(D_2) G(D_2)^2}
  \sum_{D_1\mid (P/D_2)}
\frac{H(D_1)^2 }{D_1 G(D_1)^2} \,.
\end{align*}
If we denote the innermost sum  by $I(P/D_2)$, then
\begin{align}\label{I formula}
 I(P/D_2)=\prod_{p\mid (P/D_2)}\left(1  + \frac{H(p)^2 }{p G(p)^2}\right)
 \,, \end{align}
and we find that
 \begin{align*}
 \T_2
& \ll    G(P)^2 I(P)\frac{\phi(P) }{P} \log^2 T \log^2 Y
\sum_{D_2\mid P_2}    \frac{1}{D_2  \phi(D_2) G(D_2)^2 I(D_2)} \\
&\ll    \epsilon_1^2 G(P)^2 I(P)\frac{\phi(P) }{P} \log^4 T
\prod_{p\mid P_2}   \left(1 +  \frac{1}{p  \phi(p) G(p)^2 I(p)} \right)
 \end{align*}

Now, by the definitions of $G$, $H$, and $I$ in \eqref{G formula}, \eqref{H formula},
and \eqref{I formula}, we have
$G(p)=  1-\frac{2}{p} + O(\frac{1}{p^2}) $,
$H(p)=  1+O(\frac{1 }{p}) $, and
$I(p)=  1+ \frac{1}{p} \left( (1+O(\frac{1 }{p}))/
(1-\frac{2}{p} + O(\frac{1}{p^2})) \right)^2 =  1+ \frac{1}{p}+O(\frac{1}{p^2}) $.
From these estimates it is clear  that the product over $p$ dividing $P_2$ here is
$\prod_{\sqrt{X}<p\leq X}\left( 1+ O(1/p^2) \right) \ll 1$.
Thus
\begin{align*}
  \T_2 & \ll   \epsilon_1^2 G(P)^2 I(P)\frac{\phi(P) }{P} \log^4 T \\
  & \ll  \epsilon_1^2 \log^4 T \prod_{p\mid P}\left(
  \left(1- \frac{4}{p}+ O(1/p^2) \right)
   \left(1+ O(1/p^2) \right)
  \right) \\
& \ll    \epsilon_1^2 \log^4 T \prod_{p\mid P}
  \left(1- \frac{1}{p}\right)^{4}
  \ll {\epsilon_1}^2 \left( \frac{\log T}{\log X} \right)^4 \,.
\end{align*}

The treatment of $\T_3$ is almost identical and
leads to the same bound.
Thus, combining our estimates for  $\T_1$ (see  \eqref{T_1 estimate}),
$\T_2$, and $\T_3$  with
  \eqref{Jose's Formula 2}, and noting that we may take $\epsilon_1>0$
 as small as we like, we obtain \eqref{k=2 splitting}. This
 completes the proof of the case $k=2$ of
 Theorem~\ref{thm:first case of splitting} and thus, also  the proof of the theorem.

\appendix

\section{Graphs}

\begin{figure}[p]
\centering
\includegraphics[width=3.9in]{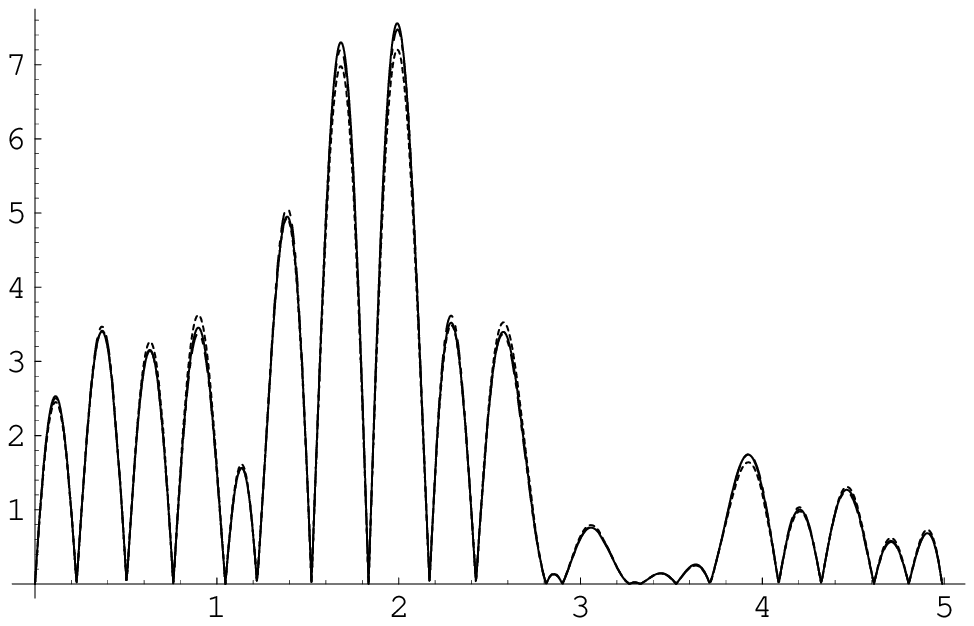}
\caption{Graph of $|\zeta(\tfrac{1}{2}+\i (x+t_0))|$ (solid) and
$|P_X(\frac12+\i (x+t_0))Z_X(\frac12+\i (x+t_0))|$, with
$t_0=\gamma_{10^{12}+40}$, with $X=\log t_0$ (dots) and $X=1000$
(dash-dots).}\label{fig:mod_zeta}

\includegraphics[width=3.9in]{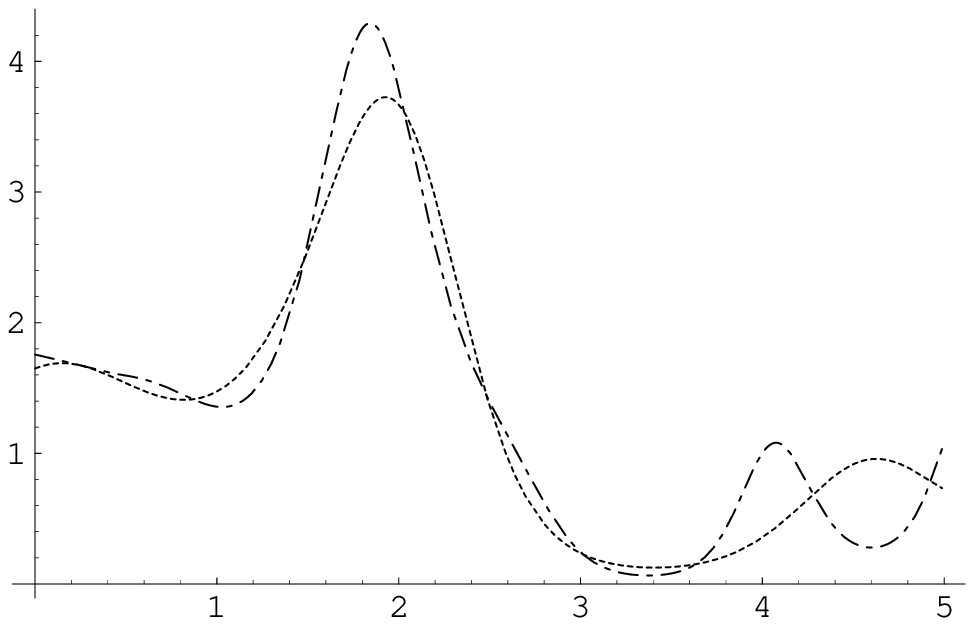}
\caption{Graph of $|P_X(\frac12+\i (x+t_0))|$, with
$t_0=\gamma_{10^{12}+40}$, with $X=\log t_0$ (dots) and $X=1000$
(dash-dots).}\label{fig:P}

\includegraphics[width=3.9in]{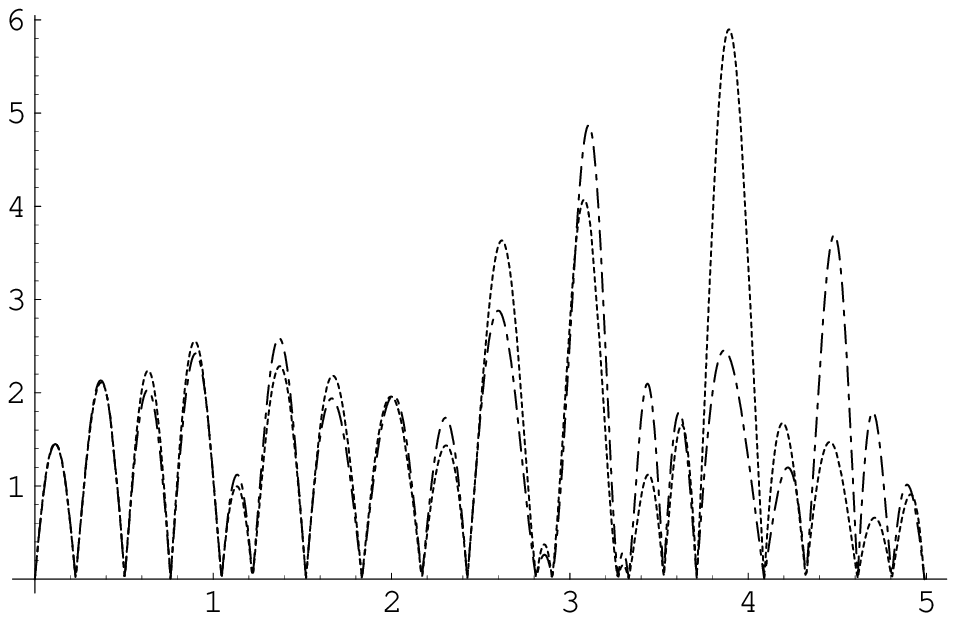}
\caption{Graph of $|Z_X(\frac12+\i (x+t_0))|$, with
$t_0=\gamma_{10^{12}+40}$, with $X=\log t_0$ (dots) and $X=1000$
(dash-dots).}\label{fig:Z}
\end{figure}

To illustrate Theorem~\ref{thm:zeta as product}, in
Figures~\ref{fig:mod_zeta}--\ref{fig:Z} we have plotted $|Z_X(\tfrac12+\i
t)|$ and $|P_X(\tfrac12+\i t)|$ for $t$ near the $10^{12}$th zero
for two values of $X$,  and have compared their product with the Riemann
zeta function. The values of $X$   used are $X= 26.31 \approx \log
\gamma_{10^{12}}$ and $X=1000$. Though the functions $P_X$ and
$Z_X$ depend upon $X$, when multiplied together the $X$ dependence
mostly cancels out, and we have an accurate pointwise
approximation to the zeta function.
The actual functions plotted are
\begin{equation*}
\left|P_X\left(\tfrac12+\i(x+t_0)\right)\right| =
\exp\left(\sum_{n\leq X} \frac{\Lambda(n) \cos((x+t_0)\log
n)}{\log n \sqrt{n}} \right)
\end{equation*}
and
\begin{equation*}
\left|Z_X\left(\tfrac12+\i(x+t_0)\right)\right|  =
\prod_{n=N+1}^{N+100} \exp\left(\Ci(|x+t_0-\gamma_n|\log X)\right)\,,
\end{equation*}
where $t_0 = \gamma_{10^{12}+40}$.
The  values of the zeros of the zeta function came from
Andrew Odlyzko's tables \cite{Odl}. The functions were plotted
for $x$ between $0$ and $5$, a range   covering the zeros
between $\gamma_{10^{12}+40}$ and $\gamma_{10^{12}+60}$. Note that
the function $Z_X$ we have plotted is an unsmoothed, truncated form of
the function
$Z_X$ that appears in Theorem~\ref{thm:zeta as product}.

\section*{Acknowledgment}
We are very grateful to Kannan Soundararajan for several comments
which strengthened some of the theorems appearing in an earlier
version of this paper.  Work of the first author was supported by
NSF grant DMS 0201457 and by an NSF Focused Research Group grant
(DMS 0244660). Work of the second author was partially supported
by EPSRC grant N09176 and an NSF Focused Research Group grant (DMS
0244660). The third author is supported by an EPSRC Senior
Research Fellowship. All three authors wish to thank the American
Institute of Mathematics and the Isaac Newton Institute for their
hospitality during the course of this work.

\end{document}